\documentclass[12pt,twoside,reqno,psamsfonts]{amsart}
\usepackage{amssymb}
\makeindex
\def\R{{\mathbb R}}

\def\Z{{\mathbb Z}}

\def\squareforqed{\hbox{\rlap{$\sqcap$}$\sqcup$}}
\def\qed{\ifmmode\squareforqed\else{\unskip\nobreak\hfil
\penalty50\hskip1em\null\nobreak\hfil\squareforqed
\parfillskip=0pt\finalhyphendemerits=0\endgraf}\fi}
\newtheorem{thm}{Theorem}[section]
\newtheorem{cor}[thm]{Corollary}
\newtheorem{lem}[thm]{Lemma}
\newtheorem{prop}[thm]{Proposition}
\DeclareMathOperator*{\ra}{\rightarrow}

\begin{document}
\title{Asymptotic boundary forms for tight Gabor frames 
and lattice localization domains}
\author{ H.G. Feichtinger}
\address{Faculty of Mathematics, University Vienna, Oskar-Morgenstern-Platz 1, 1090 Wien, AUSTRIA}
\author{K. Nowak}
\address{Department of Computer Science,
Drexel University, 3141 Chestnut Street, Philadelphia, PA 19104, USA}
\author{M. Pap}
\address{Faculty of Sciences, University of P\'ecs, 7634 P\'ecs, Ifj\'us\'ag \'ut 6, HUNGARY}
\begin{abstract}
We consider Gabor localization operators $G_{\phi,\Omega}$ defined by two
parameters, the generating function $\phi$ of a tight Gabor frame $\{\phi_\lambda\}_{\lambda \in \Lambda}$, parametrized by the elements of a given lattice $\Lambda \subset \R^2$, i.e. a discrete cocompact subgroup of $\R^2$, and a lattice localization domain $\Omega \subset \R^2$ with its boundary consisting of line segments connecting points of $\Lambda$. We find an explicit formula for the boundary form $BF(\phi,\Omega)=\text{A}_\Lambda
\lim_{R\rightarrow \infty}\frac{PF(G_{\phi,R\Omega})}{R}$, the normalized limit of the projection functional 
$PF(G_{\phi,\Omega})=\sum_{i=0}^{\infty}\lambda_i(G_{\phi,\Omega})(1-\lambda_i(G_{\phi,\Omega}))$, where $\lambda_i(G_{\phi,\Omega})$ are the eigenvalues of the localization operators $G_{\phi,\Omega}$ applied to dilated domains $R\Omega$, $R$ is an integer and $\text{A}_\Lambda$ is the area of the fundamental domain of the lattice $\Lambda$. 
Although the lattice $\Lambda$ is also a parameter of the localization operator we assume that it is fixed and we do not list it explicitly in our notation. 
The boundary form expresses quantitatively the asymptotic interactions between the generating function $\phi$ of a tight Gabor frame and the oriented boundary $\partial \Omega$ of a lattice localization domain from the point of view of the projection functional, which measures to what degree a given trace class operator fails to be an orthogonal projection. It provides an evaluation framework for finding the best asymptotic matching between pairs consisting of generating functions $\phi$ and lattice domains $\Omega$. It takes into account
directional information involving outer normal vectors of the linear segments 
constituting the boundary of $\Omega$ and the weighted sums over the corresponding half spaces of the absolute value squares of the reproducing kernels obtained out of $\phi$. In the context of tight Gabor frames and 
Gabor localization domains placing an upper bound on the square of the $L^2$ norm of the generating function $\phi$ divided by the area of the fundamental domain of lattice $\Lambda$ corresponds to
keeping the relative redundancy of the frame bounded. Keeping the area of 
the localization domain $\Omega$ bounded above corresponds to controlling the relative dimensionality of the localization problem. 
\end{abstract}
\keywords{Toeplitz operators, phase space localization, tight Gabor frames, semi-classical limit}
\subjclass{47B35,42C15,81Q20}

\maketitle
\markboth{Asymptotic Boundary Forms}{H.G. Feichtinger, K. Nowak, M. Pap}

\section{Main Results and Their Context}

We start by formulating the main results of this paper. Let $\Lambda\subset \R^2$ be a lattice,  i.e. a discrete cocompact subgroup of $\R^2$, satisfying condition $\text{A}_\Lambda<1$, where $\text{A}_\Lambda$ denotes the 
area of the fundamental domain of $\Lambda$. For a function $\phi \in L^2(\R)$ and $\lambda=(q,p)$ we define $\phi_\lambda(x)=e^{2\pi i  px}\phi(x-q)$. 
Let the system $\left\{\phi_\lambda \right\}_{\lambda \in \Lambda}$, where 
$\phi \in L^2(\R)$ is a fixed function, be a {\it tight Gabor frame}, i.e. we assume that for every $f \in L^2(\R)$ 
\begin{equation}
f=\sum_{\lambda \in \Lambda} \langle f,\phi_\lambda\rangle \phi_\lambda,
\label{gabor_reproducing_formula_discrete}\end{equation}
where the convergence of the sum is understood in the unconditional norm sense. Function $\phi$ is called the {\it generating function} of a tight Gabor frame $\left\{\phi_\lambda \right\}_{\lambda \in \Lambda}$. 
An operator $G_{\phi,b}$, acting on $L^2(\R)$, and of the form 
\begin{equation}
G_{\phi,b}f=\sum_{\lambda \in \Lambda} b(\lambda) \langle  f,\phi_\lambda\rangle \phi_\lambda
\label{gabor_multiplier}\end{equation}
is called a {\it Gabor multiplier of localization type}, if its {\it symbol} $b$, defined on the lattice $\Lambda$, is non-negative and summable. It can be easily verified that Gabor multipliers of localization type are non-negative, trace class, and that
$$
\text{tr}(G_{\phi,b})=||\phi||_{L^2(\R)}^2\sum_{\lambda \in \Lambda}b(\lambda),
$$
and the operator norm of $G_{\phi,b}$ satisfies
$$
||G_{\phi,b}||\le ||b||_{l^\infty(\Lambda)}.
$$
The {\it projection functional} $PF$ is defined on positive definite, trace class operators $T$, with their operator norm bounded above by $1$, via the formula 
\begin{equation}
PF(T)=\sum_{i=0}^{\infty}\lambda_i(T)
(1-\lambda_i(T)),
\label{def_pf}\end{equation}
where $\{\lambda_i(T)\}_{i=0}^\infty$ are the eigenvalues of $T$.
The projection functional measures the extend by which the operator $T$ 
fails to be an orthogonal projection. It takes 
non-negative values and it vanishes on the space consisting of finite dimensional orthogonal projections.
A collection of line segments $l_i\subset \R^2$, $i=1,2,...,C$ each of them starting and
ending at a point of $\Lambda$ is called a {\it $\Lambda$ cycle} if the union 
$l_1\cup l_2\cup ... \cup l_C$ forms a closed continuous line without self 
intersections. A bounded, connected and closed subset $\Omega \subset \R^2$ is called a {\it $\Lambda$ domain} if its boundary consists of a finite family of $\Lambda$ cycles $C_i$, $i=1,2,...,B$ satisfying condition $\text{dist}(C_i,C_j)>0$ for $i\ne j$. 
 We say that a function $\phi\in L^2(\R)$ satisfies {\it condition $\Phi$} if
\begin{description}
\item [$\Phi$]  \hskip4.0truecm
$\displaystyle{
\sum_{\lambda \in \Lambda}|\langle \phi,\phi_\lambda\rangle|^2
|\lambda |<\infty}$.
\label{condition_Phi} \end{description}
For a $\Lambda$ domain $\Omega \subset \R^2$ and a generating function 
$\phi$ of a tight Gabor frame $\{\phi_\lambda\}_{\lambda\in \Lambda}$ satisfying condition $\Phi$ the {\it boundary form} $\text{BF}(\phi,\Omega)$ is defined by the formula
\begin{equation}
BF(\phi,\Omega)=
\sum_{i=1}^N \text{length}(l_i)\sum_{\lambda \in U_i} \text{dist}(\lambda,P_i)|\langle \phi,\phi_\lambda\rangle|^2,
\label{definition_bf}\end{equation}
where $l_i$, $i=1,...,N$ are the line segments constituting the boundary of $\Omega$, $n_i$ is the unit vector orthogonal to $l_i$ directed outside $\Omega$, $P_i=\{w\in \R^2 \,|\,w\cdot n_i=0 \}$, $U_i=\{\lambda\in \Lambda \,|\,
\lambda\cdot n_i \ge 0 \}$.

The projective metaplectic representation defined on $SL(2,\R)$ provides the most natural way to deal with linear changes of coordinates of the time-frequency plane $\R^2$. Its definition and its basic properties are reviewed in Section 4. It came to us as a surprise, that the boundary form $BF$, representing the limit value of the projection functional $PF$ with respect to dilation factors of the localization domain tending to infinity, is invariant with respect to the action of the projective metaplectic representation. Many geometric and numerical studies of lattices in $\R^2$ distinguish the hexagonal lattice  as being special, but this is not the case in the context of the projection functional $PF$ and its limit value $BF$. 

In our first result we describe the invariance properties of $PF$ and $BF$ with respect to the projective metaplectic representation. Although the invariance of $PF$ is standard we include it as well for the sake of completeness.     
\begin{thm}\label{BF_and_PF_invariance} 
Let $\mu$ be the projective metaplectic representation of $SL(2,\R)$ acting on
$L^2(\R). $Let $\phi$ be a generating function of a tight Gabor frame $\{\phi_\lambda\}_{\lambda \in \Lambda}$ satisfying condition $\Phi$ with respect to $\Lambda$. Then for any $A\in SL(2,\R)$ $\Gamma =A(\Lambda)$ is a lattice, $\{(\mu(A)\phi)_\gamma\}_{\gamma \in \Gamma}$ is a tight Gabor frame with the generating function $\mu(A)\phi$ satisfying condition $\Phi$
with respect to $\Gamma$, and
\newline
(i) for any Gabor multiplier $G_{\phi,b}$ of localization type with the symbol bounded above by $1$, $G_{\mu(A)\phi,b\circ A^{-1}}$ is a Gabor multiplier of localization type with the symbol bounded above by $1$, and
\begin{equation}
PF\left(G_{\phi,b}\right)=PF\left(G_{\mu(A)\phi,b\circ A^{-1}}\right),
\label{pf_tranformation_property_gabor_multiplier}
\end{equation} 
\newline
(ii) for any $\Omega$ a $\Lambda$ domain, $A(\Omega)$ is a $\Gamma$ domain, and 
\begin{equation}
BF(\phi,\Omega)=BF(\mu(A)\phi,A(\Omega)).
\label{transformation_property_boundary_form}
\end{equation}
\end{thm}

Our second result is the principal result of the current paper. It describes the limit behavior of the projection functional $PF$ applied to a Gabor multiplier with the symbol of the form $b_R=\chi_{R\Omega}$, i.e. the characteristic function of a dilated lattice domain $\Omega$ by a factor $R\in \Z$, with $R\rightarrow \infty$. Gabor multipliers of this special form are called {\it Gabor localization operators} and they are denoted as $G_{\phi,R\Omega}$.  
We will see later on that the dilation factor $R$ can in fact take any real values. 
\begin{thm}\label{BF_as_PF_limit} 
Let $\phi$ be a generating function of a tight Gabor frame parametrized by lattice 
$\Lambda$ and satisfying condition $\Phi$, and $\Omega$ a $\Lambda$ domain 
contained in $\R^2$. Then
\begin{equation}
\lim_{R \rightarrow \infty}\frac{PF(G_{\phi, R\Omega})}{R}=
\frac{1}{\text{A}_\Lambda}
BF(\phi, \Omega).
\label{bf_as_pf_limit}\end{equation}
\end{thm}
\noindent
Formula (\ref{bf_as_pf_limit}) expresses the limit behavior of the projection functional $\text{PF}$ in terms of the boundary form $\text{BF}$. It provides a very explicit, quantitative way of describing the interactions between the boundary $\partial \Omega$ and the reproducing kernel obtained out of the generating function $\phi$.

The invariance of the projection functional $\text{PF}$ and the boundary form $\text{BF}$ with respect to the action of the projective metaplectic representation expressed in Theorem \ref{BF_and_PF_invariance}, taken together with the 
limit result of Theorem \ref{BF_as_PF_limit}, has its important consequences. We can conclude that no lattice $\Lambda$ is distinguished, neither from the point of view the value of the limit of $PF(G_{\phi,\R\Omega})$ as $R\rightarrow \infty$, nor from the point of view of its rate of convergence.
\begin{cor}\label{no_best_lattice}
For any lattice $\Lambda \subset \R^2$ satisfying condition $\text{A}_\Lambda <1$, any generating function $\phi$ of a tight Gabor frame $\{ \phi_\lambda \}_{\lambda\in \Lambda}$, any $\Lambda$ lattice domain $\Omega$, and any $a,b>0$ satisfying $ab=\text{A}_\Lambda$, there are a generating function $\breve{\phi}$ of a tight Gabor frame $\{ \breve{\phi}_\lambda \}_{\lambda\in a\Z \times b\Z}$ and a $a\Z \times b\Z$ lattice domain $\breve{\Omega}$, satisfying $||\phi||_{L^2(\R)}=||\breve{\phi}||_{L^2(\R)}$, $\text{Area}( {\Omega})=\text{Area}( \breve{\Omega})$, $PF(G_{\phi,R\Omega})=
PF(G_{\breve{\phi},R\breve{\Omega}})$, for all $R>0$, and also 
$BF(\phi, \Omega)=BF(\breve{\phi}, \breve{\Omega})$.
The rates of convergence of $\frac{PF(G_{\phi,R\Omega})}{R}$ to 
$\frac{1}{A_\Lambda}BF(\phi,\Omega)$, and $\frac{PF(G_{\breve{\phi},R\breve{\Omega}})}{R}$ to 
$\frac{1}{ab}BF(\breve{\phi},\breve{\Omega})$ are the same.
\end{cor}

Many important mathematical theories started in the one dimensional setup, where a multitude of additional tools is available, and then through various stages of evolution came up into their full form in any finite dimension. This was the case of the representation theory of semi-simple Lie groups, which started with the listing of all irreducible representations of $SL(2,\R)$. Classical time-frequency localization operators allow explicit diagonalization in one dimension, but not in higher dimensions, yet it was possible to transfer a large portion of one dimensional results to higher dimensions.  In one dimension, for several classes of potentials, Schr\"odinger operators can be treated via explicit formulae, yet many results that follow the 
guidelines of the one dimensional setup are also true in higher dimensions, although explicit formulae are no longer available. The development of the theory of Gabor frames follows a similar path. As this is the case of many relatively recent theories, some of its branches are still at an early, one dimensional stage, e.g. the treatment of Gabor frames with maximal lattice parameter frame set via totally positive functions of finite type done recently by Gr\"ochenig and St\"ockler in \cite{GroSto}. In many instances the distinction between one and higher dimensions is related to the differences between the theories of one and several complex variables. 
The development of the phase space theory of reproducing formulae shows many similarities with the development of the theory of Gabor frames. 
The semidirect product $\R^{2n}\rtimes Sp(n,\R)$, where $Sp(n,R)$ is the symplectic group consisting of $2n \times 2n$ invertible matrices preserving the symplectic form, together with the extended projective metaplectic representation defined on it, constitute the group of affine transformations of the phase space, acting geometrically on $\R^{2n}$, and analytically on  $L^2(\R^n)$. The affine transformations of the phase space  provide a natural framework for the constructions of reproducing formulae. All reproducing formulae coming out of connected Lie subgroups of  $\R^{2n}\rtimes Sp(n,\R)$ were characterized by De Mari and Nowak for $n=1$ in \cite{DMNo}, but no analogous results are known in higher dimensions. In dimension one $Sp(1,\R)=SL(2,\R)$ and it is possible to obtain the list of all connected Lie subgroups of $\R^2\rtimes SL(2,\R)$ out of the classical structure 
results describing the orbits of inner automorphisms of $SL(2,\R)$. No analogous lists of subgroups are available in higher dimensions. 
In the current paper we deal with tight Gabor frames parametrized by a lattice $\Lambda \subset \R^2$. Any 
tight Gabor frame parameterized by a lattice $\Lambda \subset \R^2$ can be transferred into a tight Gabor frame parameterized by a separable lattice $a\Z\times b\Z$, $a,b>0$, by an appropriate $SL(2,\R)$ linear transformation of the time-frequency plane. We treat the case of a separable lattice via an explicit computation and then we transfer the formula we obtain to an arbitrary lattice with the  help of the projective metaplectic representation. This approach does not generalize to higher dimensions. 

The principal results of the current paper deal with 
discrete one-dimensional setup of tight Gabor frames, therefore for the sake of consistency we formulate  
definitions and reference results only in one dimension. Operators of composition of convolution with $g$ followed by a pointwise multiplication by $f$, where both functions $f$ and $g$ are defined on the real line $\R$, have integral kernels of the form
\begin{equation}
K(x,y)=f(x)g(x-y).
\label{convolution_product_kernel}\end{equation}
They are commonly called convolution-product operators. 
Historically three dimensional convolution-product operators played an important role in the study of Schr\"odinger operators. The Birman-Schwinger principle 
allows a transition from a  Schr\"odinger operator
to a convolution-product operator. Out of that transition it was possible to obtain sharp estimates  for the number of bound states of the Schr\"odinger operator. Classical time-frequency localization operators are one dimensional convolution-product operators with $f$ and $\check{g}$ characteristic functions of intervals, where $\check{g}$ is the inverse Fourier transform of $g$. Their spectral properties were carefully studied 
many years ago by  Landau, Pollak, Slepian and Widom (see \cite{Dau}, \cite{LanWid}, \cite{Lan}, and the references provided within), yet till now these classical  results  bring important ingredients for both theoretical and
applied components occurring in many recent developments. There were also extensive studies of convolution-product operators with one of the function parameters $f$ or $g$ fixed and of prescribed potential type (see the book by Mazya \cite{MazSha}), but the dependence on both parameters, and the mutual interaction
between $f$ and $g$ in the general case seems to be a difficult problem that is to a large degree still open. 

We are interested in operators with integral kernels of the form  (\ref{convolution_product_kernel}), where translations constituting the convolution with $g$ are  extended to combined actions of translations and modulations of the Schr\"odinger representation, and applied to a function $\phi$. As a consequence of this extension we need to introduce two arguments of the multiplier $f$, one corresponding to translations, the other to modulations, and the extended integral kernel (\ref{convolution_product_kernel}) becomes
\begin{equation}
K((q,p),y)=f(q,p)\overline{\phi_{q,p}(y)},
\label{square_root_gabor_toeplitz_kernel}\end{equation}
where $\phi_{q,p}(x)=e^{2\pi i px}\phi(x-q)$. We add complex conjugate 
over $\phi_{q,p}$, because
we want our operators to be exactly square roots of Gabor-Toeplitz operators. The Gabor reproducing formula has the form
\begin{equation}
f=\int_{\R^2} \langle f,\phi_{q,p}\rangle \phi_{q,p}\,dq\,dp ,
\label{gabor_reproducing_formula}\end{equation}
where $f,\phi \in L^2(\R)$, $||\phi ||_{L^2(\R)}=1$, and the convergence of the 
integral is understood in the weak sense. Introducing 
into  (\ref{gabor_reproducing_formula}) a weight function $b(q,p)$, 
called a symbol we obtain a Gabor-Toeplitz operator
\begin{equation}
T_{\phi,b}f=\int_{\R^2} b(q,p)\langle f,\phi_{q,p}\rangle \phi_{q,p}\,dq\,dp.
\label{gabor_toeplitz_operator}\end{equation}
Gabor-Toeplitz operators generalize Fock space Toeplitz operators. The Bargmann transform provides their mutual unitary equivalence, in case the normalized Gaussian is chosen for the generating function. Books by Folland \cite{Fol} and Zhu \cite{Zhu} are very good references on the subject. In the field of phase space analysis Gabor-Toeplitz operators were introduced by Ingrid Daubechies. Her book \cite{Dau} provides a
comprehensive account of the background and the initial results. 
It is possible to study Gabor-Toeplitz operators in a very general context, with various classes of symbols and acting on a wide range of function spaces. It is convenient to discuss Gabor multipliers and Gabor-Toeplitz operators in parallel.
In this paper we restrict attention to Gabor-Toeplitz operators of localization type acting on $L^2(\R)$, i.e. we assume that the symbol $b$ is non-negative, bounded and integrable with respect to 
the Lebesgue measure on $\R^2$. It is straightforward to verify that under these 
assumptions $T_{\phi,b}$ is non-negative, trace class,
$$
\text{tr}(T_{\phi,b})=\int_{\R^2}b(q,p)\,dq\,dp,
$$
the operator norm of $T_{\phi,b}$ satisfies
$$
||T_{\phi,b}||\le ||b||_{L^\infty(\R^2)},
$$
and that if $|f(q,p)|^2=b(q,p)$, then the composition of the operator defined by the integral kernel   (\ref{square_root_gabor_toeplitz_kernel}) with its
conjugate equals $T_{\phi,b}$ defined in (\ref{gabor_toeplitz_operator}). Kernels of the form (\ref{square_root_gabor_toeplitz_kernel}) represent square roots of Gabor-Toeplitz operators (\ref{gabor_toeplitz_operator}). They provide a link between Gabor-Toeplitz operators and generalized convolution-product operators, where translations are substituted by actions of unitary representations.

Asymptotic properties, as $R\rightarrow \infty$, of the symbolic calculus of Gabor-Toeplitz operators 
$T_{\phi,b_R}$, where $b$ is integrable, $0 \le b(q,p) \le 1$, and 
$b_R(q,p)=b(\frac{q}{R},\frac{p}{R})$ is the $L^\infty$ normalized dilation 
of $b$, were studied in \cite{FeiNo1}. 
For  $h$ a continuous function defined on the closed interval $[0,1]$ the operator
$h(T_{\phi,b_R})$ is defined via the spectral decomposition of $T_{\phi,b_R}$.
A Szeg\"o type formula for operators of the form $h(T_{\phi,b_R})$ was obtained
in \cite{FeiNo1}, showing that
\begin{equation}
\lim_{R\rightarrow \infty}\frac{\text{tr}(T_{\phi,b_R}h(T_{\phi,b_R}))}{R^2}
=\int_{\R^2}b(q,p)h(b(q,p))\,dq\,dp.
\label{szego_formula_gabor_toeplitz} \end{equation}
We need to multiply $h(T_{\phi,b_R})$ by $T_{\phi,b_R}$ in formula (\ref{szego_formula_gabor_toeplitz}) in order to ensure
that the operator is trace class. The
Bohr correspondence principle  was concluded as a consequence of  (\ref{szego_formula_gabor_toeplitz}). It asserts that in the normalized limit with factor $R^2$ in the denominator both the distribution of the eigenvalues of $T_{\phi,b_R}$
\begin{equation}
N(\delta,R)=\left| \left\{i\,|\,\lambda_i(T_{\phi,b_R}) > \delta \right\}\right|, 
\,0<\delta <1,
\label{eigenvalue_distribution_with_dilations}\end{equation}
and the size of their plunge region 
\begin{equation}
M(\delta_1,\delta_2,R)=\left| \left\{i\,|\,\delta_1 < \lambda_i(T_{\phi,b_R}) < \delta_2 \right\}\right|, \,0<\delta_1 <\delta_2 <1,
\label{eigenvalue_plunge_region_with_dilations}\end{equation}
are expressed directly via the corresponding quantities of the symbol function,
the distribution of $b$
\begin{equation}
\left| \left\{(q,p)\,|\,b(q,p) > \delta \right\}\right|, 
\label{distribution_symbol}\end{equation}
and the Lebesgue measure of its plunge region
\begin{equation}
\left| \left\{(q,p)\,|\, \delta_1 <b(q,p) <\delta_2 \right\}\right|.
\label{plunge_region_symbol}\end{equation}
In the passage to the normalized limit it is necessary to assume that the level sets
$\{(q,p)\,|\,b(q,p)=\tau\}$, $\tau=\delta, \delta_1, \delta_2$ , have Lebesgue measure $0$.
We were able to deduce that asymptotically the best localization properties, i.e. 
asymptotically there are no eigenvalues in the open interval $(0,1)$, occur for symbols being characteristic functions of measurable sets $\Omega\subset \R^2$. 

In the next step operators $T_{\phi,\Omega}$, with symbols $b=\chi_\Omega$, 
were studied directly without the asymptotic limit. 
Two-sided estimates of the size of the eigenvalue plunge region
\begin{equation}
M(\delta_1,\delta_2,\phi,\Omega)=\left| \left\{i\,|\,\delta_1 < \lambda_i(T_{\phi,\Omega}) < \delta_2 \right\}\right|
\label{plunge_region_domains_continuous}\end{equation}
expressed in terms of the area of a strip of fixed size $R$ around the boundary 
$\partial \Omega$ and uniform with respect to generating functions $\phi$ and localization domains $\Omega$, of the form
\begin{equation}
c_1\left| (\partial \Omega)^R\right| \le M(\delta_1,\delta_2,\phi,\Omega)
\le c_2\left| (\partial \Omega)^R\right|,
\label{plunge_region_estimates_domains_continuous}\end{equation}
where $(\partial \Omega)^R=\{(q,p)\,|\,\text{dist}((q,p),\partial \Omega)| < R\}$, were obtained in \cite{FeiNo2}. In order to get two-sided uniform 
positive constants $c_1,c_2$
it was necessary to assume uniform decay and non-degeneracy of the reproducing kernels obtained out of generating functions $\phi$ and uniform access 
to localization domains $\Omega$ and their complements $\Omega^c$ from the points near 
their boundaries $\partial\Omega$. 
The size $R$ of the strip  $(\partial \Omega)^R$ around $\partial \Omega$ is one of the uniform parameters controlling mutual interactions between the generating functions $\phi$
and domains of localization $\Omega$. 

The next step in the study of mutual interactions between generating functions
and domains of localization from the point of view of estimating the size 
of the eigenvalue plunge region was accomplished in \cite{No1}.
The eigenvalues of Gabor-Toeplitz localization operators 
$\lambda_i(T_{\phi,\Omega,})$ 
satisfy estimates  $0\le\lambda_i(T_{\phi,\Omega,})\le 1$.
The projection functional 
$PF(T_{\phi,\Omega})=\sum_{i=0}^\infty
\lambda_i(T_{\phi,\Omega,})(1-\lambda_i(T_{\phi,\Omega,}))$
provides an exact, quantitative way of measuring the size of the eigenvalue plunge region. It has the same form as in (\ref{def_pf}) in the discrete setup.
We assume that $\Omega$ is a bounded domain with $C^1$ boundary.
Symbol $\sigma$ denotes the arc length defined on $\partial \Omega$. Function
$n$ is the Gauss map, i.e. $n(r,s)$ is the unit normal vector 
at $(r,s)\in \partial \Omega$ directed outside $\Omega$. By 
$P_v=\{w\in \R^2\,|\, w \cdot v=0\}$ we denote the linear subspace of $\R^2$ consisting
of vectors orthogonal to $v$, and $U_v=\{w\in \R^2\,|\,w \cdot v \ge 0\}$ is the half space
inside $\R^2$ with $v$ being the inner normal vector. The boundary form 
is defined as
\begin{equation}
BF(\phi,\Omega)=
\int_{\partial \Omega}\int_{U_{n(r,s)}}|\langle \phi,\phi_{q,p}\rangle |^2
\text{dist}((q,p),
P_{n(r,s)})dq\,dp\,d\sigma (r,s).
\label{def_bf}\end{equation}
Under an appropriate integrability condition imposed on $\phi$, condition 
$\Phi$ of the current paper is its adaptation to the discrete setup, that makes formula (\ref{def_bf}) 
well defined it has been shown in \cite{No1} that
\begin{equation}
\lim_{R\rightarrow \infty}\frac{PF(T_{\phi,R\Omega,})}{R}=
BF(\phi, \Omega).
\label{def_af}\end{equation}
Formula ($\ref{def_af}$) expresses the limit behavior of the projection functional 
in terms of the boundary form. It describes quantitatively the interactions between the boundary $\partial \Omega$ and the 
reproducing kernel obtained out of the generating function $\phi$. Boundary form (\ref{definition_bf}) is an adaptation of (\ref{def_bf}) to the discrete setup.

Let us fix the area $A$ and a generating function $\phi$ for which the boundary form (\ref{def_bf}) is well defined,  
and let us ask for what localization domains $\Omega$ of area $A$ the boundary form $BF(\phi,\Omega)$ takes the smallest possible value.
Let us introduce surface tension $M_\phi$, 
defined on the unit sphere of the time-frequency plane $\R^{2}$, 
given by the formula 
\begin{equation}
M_\phi (v)=
\int_{U_v}|\langle \phi,\phi_{q,p}\rangle |^2\text{dist}((q,p),
P_v)dq\,dp.
\label{def_st}\end{equation}
Function $M_\phi$, defined in (\ref{def_st}), determines the Wulff shape  
\begin{equation}
K_\phi=\bigcap_{v\in S^{1}}\{w\in \R^{2}\,|\,w \cdot v \le M_\phi(v)\}.
\label{def_ws}\end{equation}
For the given generating function $\phi$ let us construct $K_\phi$, defined in 
(\ref{def_ws}), and 
let us scale it with a constant $c$, so that the area of $cK_\phi$ is $A$. 
It occurs, that, up to translation, the domain $cK_\phi$ is the only domain
among all domains $\Omega$ with finite perimeter and area $A$ for which $BF(\Omega,\phi)$ is minimal.  
The shape $cK_\phi$ is the optimal domain of localization.

The principal results of the current paper deal with the discrete setup 
of tight Gabor frames defined on $L^2(\R)$. On the level of convolution-product operators  with kernels of the form (\ref{convolution_product_kernel}) the
transition to the discrete setup translates to making the range variables discrete and keeping the domain variables continuous. The multiplication parameter $f$ is now defined on the group of integers, but the convolution parameter $g$ is still defined on the real line $\R$. After the modification the integral kernel (\ref{convolution_product_kernel}) becomes
\begin{equation}
K(m,y)=f(m)g(m-y).
\label{convolution_product_kernel_discrete}\end{equation}
In the discrete setup the Gabor reproducing formula (\ref{gabor_reproducing_formula}) is 
substituted by a tight Gabor frame expansion (\ref{gabor_reproducing_formula_discrete})
$$
f=\sum_{\lambda \in \Lambda} \langle f,\phi_\lambda\rangle \phi_\lambda,
$$
and the Gabor-Toeplitz operator
(\ref{gabor_toeplitz_operator}) by a Gabor multiplier (\ref{gabor_multiplier})
$$
G_{\phi,b}f=\sum_{\lambda \in \Lambda} b(\lambda) \langle  f,\phi_\lambda\rangle \phi_\lambda.
$$
Gabor multipliers were introduced as a phase space analysis tool parallel to Gabor
expansions. Both Gabor-Toeplitz operators and Gabor multipliers are currently very actively investigated, mostly from the point of view of their usage in phase space analysis. They were applied as phase space partitioning operators in \cite{DorGro} , leading to a new characterization of modulation spaces, as isomorphism maps,
in \cite{GroTof1}, \cite{GroTof2}, in the context of modulation spaces and weighted Bargmann-Fock spaces, in \cite{Dop}, \cite{Gro2}, \cite{DorTor}, \cite{CorNicRod}, \cite{CorGroNic} as approximation blocks for the representation of Hilbert-Schmidt, pseudodifferential, and Fourier integral operators. A recent survey \cite{FeiNoPa} provides an overview of the results on both Gabor-Toeplitz operators and Gabor multipliers from the point of view of phase space localization. Not all of the phenomena are the same for Gabor-Toplitz operators and Gabor multipliers. The cut-off phenomenon happens only in the continuous setup (see \cite{No2}), the Berezin transform is invertible in a stable way only in the discrete setup (see \cite{FeiNo3}).

We restrict attention to Gabor multipliers of localization type, i.e. we assume that the symbol $b$ defined on $\Lambda$ is non-negative and summable. 
The operator given by the kernel
$$
K(\lambda,y)=f(\lambda)\overline{\phi_\lambda(y)},
$$
a discrete analogue of (\ref{square_root_gabor_toeplitz_kernel}), acting from $L^2(\R)$ into $l^2(\Lambda)$,
is again a square root of $G_{\phi,b}$, provided $|f(\lambda)|^2=b(\lambda)$.

A discrete analogue of the Szeg\"o type formula (\ref{szego_formula_gabor_toeplitz}) applied to $G_{\phi,b_R}$ required us to assume that the symbol function
$b$ is defined on $\R^2$, that it is Riemann integrable, has compact support, and
that $0 \le b(q,p) \le 1$.  Also in the discrete setup the Szeg\"o type formula captures essential asymptotic properties of the symbolic 
calculus of Gabor multipliers of localization type as $R\rightarrow \infty$, where $R$
is the parameter of the $L^\infty$ normalized dilation $b_R(q,p)=b(\frac{q}{R},\frac{p}{R})$. 
For any continuous function $h$ defined on $[0,1]$ the operator $h(G_{\phi,b_R})$
is defined in terms of the spectral decomposition of $G_{\phi,b_R}$. The discrete analogue of (\ref{szego_formula_gabor_toeplitz}) was obtained in \cite{FeiNo3}.
It states that
\begin{equation}
\lim_{R\rightarrow \infty}\frac{\text{tr}(G_{\phi,b_R}h(G_{\phi,b_R}))}{R^2}
=\frac{||\phi||_{L^2(\R)}^2}{\text{A}_\Lambda}\int_{\R^2}b(q,p)h(b(q,p))\,dq\,dp.
\label{szego_formula_gabor_multipliers} \end{equation}
The motivation for both formulae (\ref{szego_formula_gabor_toeplitz}) and
(\ref{szego_formula_gabor_multipliers}) came out of 
Widom's paper \cite{Wid}. The Bohr correspondence principle, i.e. 
the normalized limits, the same as for the eigenvalue distribution of
(\ref{eigenvalue_distribution_with_dilations}),  
and the size of their plunge region of
(\ref{eigenvalue_plunge_region_with_dilations}), can be concluded for 
Gabor multipliers of localization type. In the discrete setup we need to insert the multiplicative
constant $\frac{||\phi||_{L^2(\R)}^2}{\text{A}_\Lambda}$ in front of 
(\ref{distribution_symbol}) and (\ref{plunge_region_symbol}),
the same constant as in (\ref{szego_formula_gabor_multipliers}), but otherwise 
the normalized limits are same as in the case of Gabor-Toeplitz operators of localization type.
In particular for a bounded, Riemann measurable set $\Omega$ we obtain
$$
\lim_{R\rightarrow \infty}\frac{|\{i\,|\,\lambda_i(G_{\phi,R\Omega})>
1-\epsilon\}|}{R^2}=\frac{||\phi||_{L^2(\R)}^2}{\text{A}_\Lambda}
|\Omega|,
$$
for any $\epsilon >0$, where in case $b_R=\chi_{R\Omega}$ we use the notation 
$G_{\phi,b_R}=G_{\phi,R\Omega}$. According to the principles of quantum mechanics we interpret 
the area $|\Omega|$ as the dimension of the space of functions localized in phase space to $\Omega$. The constant $\frac{||\phi||_{L^2(\R)}^2}{\text{A}_\Lambda}$ compensates the relative redundancy of the tight frame.

Two-sided uniform estimates of the size of the eigenvalue plunge region, the discrete
analog of (\ref{plunge_region_domains_continuous}), were also obtained in \cite{FeiNo3}. In the discrete context localization domains $\Omega$ are finite subsets of $\Lambda$, the strips of size $R$ around the boundary $\partial \Omega$ occurring 
in (\ref{plunge_region_estimates_domains_continuous}) are substituted by
$$
(\partial \Omega)^R=\left\{\lambda \in \Omega \,|\, \text{dist}(\lambda,\Omega^c)
<R \right\}\cup \left\{\lambda\in \Omega^c \,|\, \text{dist}(\lambda,\Omega)
<R \right\},
$$
and $|(\partial \Omega)^R|$ by the number of points of $\Lambda$ inside $(\partial \Omega)^R$. 

Out of Theorems \ref{BF_and_PF_invariance}, \ref{BF_as_PF_limit} we were able to conclude Corollary \ref{no_best_lattice}, expressing the fact that from the point of view of the projection functional $PF$ and the boundary form $BF$
no lattice $\Lambda$ is distinguished. It would be interesting to isolate those phase space phenomena that make distinction between lattices parameterizing tight Gabor frames. The book by Martinet \cite{Mar} is a comprehensive resource on lattices in Euclidean spaces. 

We would like to indicate intuitive reasons explaining the existence of the normalized limit of the projection functional applied to 
$G_{\phi,R\Omega}$. Let $F$ be the absolute 
value square of the reproducing kernel, i.e. 
$F(\nu)=|\langle \phi,\phi_\nu \rangle |^2$. 
The properties of the symbolic calculus of Gabor multipliers allow us 
to express the normalized projection functional as
\begin{equation}
\frac{1}{R} \sum_{\nu_1,\,\nu_2 \in \Lambda} \chi_{R\Omega^c}(\nu_1)F(\nu_1-\nu_2)
\chi_{R\Omega}(\nu_2).
\label{general_kernel}\end{equation}
We observe that 
$\chi_{R\Omega}(\nu)=\chi_{\Omega}(\frac{\nu}R)$ and 
$\chi_{R\Omega^c}(\nu)=\chi_{\Omega^c}(\frac{\nu}R)$, i.e.
instead of stretching the domain $R$ times and keeping the lattice fixed 
we may make the lattice $R$ times denser and keep the domain fixed. 
What counts are the relative sizes. Expression ($\ref{general_kernel}$) may be rewritten as
\begin{equation}
\frac{1}{R} \sum_{\nu_1,\,\nu_2 \in \frac{1}{R}\Lambda} 
\chi_{\Omega^c}(\nu_1)F(R(\nu_1-\nu_2))
\chi_{\Omega}(\nu_2).
\label{general_kernel_domain_fixed}\end{equation}
Form ($\ref{general_kernel_domain_fixed}$) of ($\ref{general_kernel}$) 
is more convenient for a conceptual geometric interpretation. We see that 
for large values of $R$ the main contribution comes from pairs $\nu_1,\nu_2$, 
$\nu_1 \in \Omega^c$, $\nu_2 \in \Omega$, which lie close together. Therefore, we expect that the resulting limit will be interpretable as an appropriate boundary form.

In the continuous case the Wulff shape (\ref{def_ws}) is the optimal localization domain. There is no analogue of it for lattice domains. We do not know how large is the class of localization domains for which the asymptotic boundary forms exist. We expect that in the general case the problem of existence of normalized limits of projection functionals has to be considered in parallel with 
the asymptotic properties of the counting function $|\Lambda\cap R\Omega|$
as $R\rightarrow \infty$. Harmonic analysis background on the counting 
problem $|\Lambda\cap R\Omega|$ as $R\rightarrow \infty$, together with the 
original Hlawka result, and several developments that followed afterwards are presented in Stein's book \cite{Ste}. The paper by Nowak \cite{Now} discusses more recent approaches to the topic.

Books by Christensen \cite{Chr}, Daubechies \cite{Dau}, Flandrin \cite{Fla},
Gr\"ochenig \cite{Gro1},  Wojtaszczyk \cite{Woj}, and collected volumes by 
Feichtinger, Strohmer \cite{FeiStr1}, \cite{FeiStr2} present a broad background of phase space analysis techniques needed for the theory of Gabor expansions.
Papers by Balan \cite{Bal}, Cassaza \cite{Cas}, Heil \cite{Hei} illustrate several further aspects of Gabor analysis not treated directly in this paper. 

The major part of this work was done during the special semester at the Erwin Schr\"odinger Institute in Vienna devoted to Time-Frequency Analysis. The authors express many thanks to the the Schr\"odinger Institute, for the hospitality, and for the great creative atmosphere. They also thank Monika D\"orfler for numerous insightful discussions on the subject of the current paper.

\section{Tight Gabor Frames, their Constructions and Properties}

Let $\mathcal{H}$ be a Hilbert space. A family of functions $\left\{ \phi_j\right\}_{j\in J}\subset\mathcal{H}$
is called a frame of $\mathcal{H}$ if there are constants $0<A,B<\infty$ such that for any $f\in \mathcal{H}$
\begin{equation}
A||f||^2_{\mathcal{H}}\le \sum_{j\in J}\left| \langle f,\phi_j\rangle _{\mathcal{H}}\right|^2
\le B||f||^2_{\mathcal{H}}.
\label{general_frames}\end{equation}
For any frame of $\mathcal{H}$ it is possible to choose the largest constant $A$ and the smallest 
constant $B$ for which (\ref{general_frames}) holds. These two extreme values of
$A$ and $B$ are called the lower and the upper frame bounds 
of the frame $\left\{ \phi_j\right\}_{j\in J}$. If the lower and the upper frame bounds are equal then the frame is called tight. 
Any tight frame gives rise to a discrete reproducing formula.  
We may renormalize a tight frame $\left\{ \phi_j\right\}_{j\in J}$ and obtain 
the norm equality 
 $$
||f||^2_{\mathcal{H}} = \sum_{j\in J}\left| \langle f,\phi_j\rangle _{\mathcal{H}}\right|^2,
$$ 
which in turn via the polarization identity may be interpreted as a reproducing formula
 \begin{equation}
f = \sum_{j\in J} \langle f,\phi_j\rangle _{\mathcal{H}} \phi_j,
\label{tight_gabor_frame_reproducting_formula} \end{equation}
with the convergence of the sum understood in the unconditional norm sense.
There is a canonical way of constructing a tight frame out of a frame. 
For any frame $\left\{ \phi_j\right\}_{j\in J}$ 
we define a frame operator via
\begin{equation}
Sf = \sum_{j\in J} \langle f,\phi_j\rangle _{\mathcal{H}} \phi_j.
\label{frame_operator} \end{equation}
Condition (\ref{general_frames}) guarantees that the above sum is unconditionally 
convergent for any $f\in\mathcal{H}$ and that the linear operator $S$ 
defined in (\ref{frame_operator}) is positive definite, bounded and invertible on $\mathcal{H}$. It is straightforward to verify that 
the family  $\left\{ S^{-1/2}\phi_j\right\}_{j\in J}$ is a tight frame on $\mathcal{H}$.

Gabor frames have the form $\left\{\phi_\lambda \right\}_{\lambda \in \Lambda}$,
where $\phi \in L^2(\R)$, $\phi_{q,p}(x)=e^{2\pi i  px}\phi(x-q)$, and 
$\Lambda$ is a lattice in $\R^2$. We assume that the system $\{\phi_\lambda\}_{\lambda \in \Lambda}$ 
is a tight Gabor frame normalized in such a way that the reproducing formula 
(\ref{tight_gabor_frame_reproducting_formula}) holds. If necessary we multiply
the generating function $\phi$ by an appropriate normalization constant. The principles of quantum mechanics allow us to interpret 
Gabor frames as coverings of the phase space. Each frame element $\phi_\lambda$
corresponds to a phase space block of area 1. 
Let $\text{A}_\Lambda$ denote the area of the fundamental domain of 
lattice $\Lambda$. If $\text{A}_\Lambda>1$ then the lattice $\Lambda$ is
too sparse and it is not possible to cover the phase space with blocks of 
area 1. There are no Gabor frames if $\text{A}_\Lambda>1$. 
If $\text{A}_\Lambda=1$, this is the so called critical case, 
it is possible to cover the phase space with blocks of area 1, but without any overlapping. It occurs that in this case a Gabor frame is automatically a Riesz basis,
and a normalized tight Gabor frame is an orthonormal basis. A well known Balian-Low theorem applies to the critical case. It states that the generating function of a Gabor frame has strictly limited smoothness and localization properties. 
It is possible to interpret it as the outcome of the covering property with no overlapping, i.e. a partitioning of the phase space. From the point of view of modulation spaces, which require smoothness and localization, there are no Gabor frames in the critical case $\text{A}_\Lambda=1$. If $\text{A}_\Lambda<1$ 
then the lattice $\Lambda$ is dense enough and it is possible to cover the phase
space with blocks of area 1. There is also a sufficient overlap between the blocks that allows constructions of generating functions of Gabor frames with arbitrary smoothness and localization properties as measured by the scale of modulation spaces. The book by Gr\"ochenig \cite{Gro1} is a comprehensive resource on all these facts. Further on we assume that $\text{A}_\Lambda < 1$.

Let $W_s(q,p)$, $s\ge 0$, be a weight function defined on $\R^2$ by the formula
\begin{equation}
W_s(q,p)=\left( 1+\left(|q|^2+|p|^2\right)^{1/2}\right)^s.
\label{standard_weight} \end{equation}
Let $S(\R^2)$ denote the space of Schwartz class functions and
$S'(\R^2)$ the space of tempered distributions. Let us select a non-zero 
function $g\in S(\R^2)$ and define the the modulation space $M^{p,q}_{W_s}$,
$1\le p,q \le \infty$, as the space consisting of those tempered distributions $f\in S'(\R^2)$ for which the norm
\begin{equation}
||f||^{p,q}_{W_s}=\left(\int_{\R} \left(\int_{\R}
\left| \langle f,g_{q,p}\rangle \right|^p W_s(q,p)^p\,dq
\right)^{q/p} dp \right) ^{1/q}
\label{modulation_spaces} \end{equation}
is finite. If $p=\infty$ or $q=\infty$ the integral norm is substituted by the 
essential supremum. The definition of the modulation space does not depend on the 
choice of function $g$. For different functions $g$ the corresponding norms are equivalent. If $p=q$ we write $M^p_{W_s}$ instead of $M^{p,p}_{W_s}$.

Let us recall that a function $\phi\in L^2(\R)$ satisfies condition $\Phi$ with respect to lattice $\Lambda$ if and only if
$\sum_{\lambda \in \Lambda}|\langle \phi,\phi_\lambda\rangle|^2
|\lambda |<\infty$.
Condition $\Phi$ is necessary and sufficient for the boundary form 
 $BF(\phi,\Omega)$ of formula (\ref{definition_bf})
to be well defined for any $\Lambda$ lattice domain $\Omega$. The next proposition
formulates a sufficient condition on a function $\phi \in L^2(\R)$, expressed
in terms of modulation spaces, for condition $\Phi$ to hold. 
Our sufficient condition follows by a direct application of the Young's inequality and the standard 
theory of modulation spaces presented in \cite{Gro1}. We do not include its proof.
\begin{prop} \label{sufficient_condition_for_Phi}
If $\phi\in L^2(\R)\cap M^{4/3}_{W_{1/2}}$, then $\phi$ satisfies condition $\Phi$.
\end{prop}

The first construction of tight Gabor frames in dimension $1$ was obtained by Daubechies, Grossmann and Meyer in 1986. Tight Gabor frames were called painless nonorthogonal expansions back then. The construction produced generating functions  with compact support either in position or in momentum and with arbitrary smoothness measured by the number of continuous derivatives, any $C^k$ or $C^\infty$ was possible. The book by Daubechies \cite{Dau} provides an excellent account of the initial stages of the constructions of tight Gabor frames. Then came Wexler-Raz biorthogonality relations, Walnut, Janssen representations, and Ron-Shen duality principle dealing with Gabor frames and the frame operator in any finite dimension. All these are very nicely presented in Gr\"ochenig's book \cite{Gro1}. Tight Gabor frames and in particular canonical tight Gabor frames obtained via the action of $S^{-1/2}$, where $S$ is the frame operator, are of the principal interest from the point of view of this paper. Generating functions being members of the modulation space $M^1_W$,
where $W$ is a subexponential weight, are the building blocks for all other modulation spaces. Subexponential weights are a natural generalization of the 
standard weights $W_s$ defined in (\ref{standard_weight}) (see \cite{GroLei} for the definition, provided on page 10, and for explanations of their usage in time-frequency analysis). It is important to know how to construct generating functions of tight Gabor frames, which are members of  $M^1_W$. The major result in this direction was obtained by Gr\"ochenig and Leinert in \cite{GroLei}, where they proved that, for any lattice $\Lambda\subset \R^{2n}$, the canonical tight frame operator $S^{-1/2}$ is bounded on $M^1_W$. If we have a Gabor frame with the generating function in $M^1_W$, 
then we also have a tight Gabor frame with the generating function in $M^1_W$. 
The existence of a generating function in $M^1_W$ of a Gabor frame, for any lattice $\Lambda\subset \R^{2n}$, 
with the volume of the fundamental domain satisfying $\text{Vol}_\Lambda < 1$, was established recently by Luef \cite{Lue}. Feichtinger and Kaiblinger examined, from the point of view modulation spaces, continuity properties of the canonical dual generating function, i.e. $\tilde{g}=Sg$,
where $S$ is the frame operator, with respect to the perturbations of the lattice.  They proved in \cite{FeiKai} that the set 
$\mathcal{F}_{M^1_{W}}=\{(\phi,A)\in M^1_{W} \times GL(\R^{2n})\,|\,
\{\phi_{A(k)}\}_{k\in\Z^{2n}}\text{ is a Gabor frame}\}$ is open and that the map 
$(\phi,A)\mapsto \tilde{\phi}$ is continuous from $\mathcal{F}_{M^1_{W}}$
into $M^1_{W}$. This result was recently extended to canonical tight generating functions $\breve{\phi}=S^{-1/2}(\phi)$ by Leinert and Luef \cite{Lue}. 

\section{Proofs of the Main Results}
We start by recalling the definitions of our principal objects of interest. A tight Gabor frame  $$\left\{\phi_\lambda \right\}_{\lambda \in \Lambda},$$ 
parametrized by a lattice $\Lambda\subset \R^2$, and normalized in such a way that the reproducing formula (\ref{gabor_reproducing_formula_discrete}) holds, is given. A Gabor multiplier of localization type (\ref{gabor_multiplier}) is constructed out of the reproducing system $\left\{\phi_\lambda \right\}_{\lambda \in \Lambda}$, and 
a non-negative, summable, and bounded above by $1$ symbol $b$, defined on $\Lambda$. It is defined by the formula
$$
G_{\phi,b}f=\sum_{\lambda \in \Lambda} b(\lambda) \langle  f,\phi_\lambda\rangle \phi_\lambda.
$$
The projection functional (\ref{def_pf}) applied to $G_{\phi,b}$ has the form
$$
PF(G_{\phi,b})=\sum_{i=0}^{\infty}\lambda_i(G_{\phi,b})
(1-\lambda_i(G_{\phi,b})),
$$
where $\{\lambda_i(G_{\phi,b})\}_{i=0}^\infty$ are the eigenvalues of $G_{\phi,b}$. Condition $\Phi$, the additional requirement imposed on the generating function $\phi$, 
$\sum_{\lambda \in \Lambda}|\langle \phi,\phi_\lambda\rangle|^2
|\lambda |<\infty$, is a necessary and a sufficient condition for the boundary form (\ref{definition_bf}) to be well defined for any $\Lambda$ localization domain $\Omega$,
$$
BF(\phi,\Omega)=
\sum_{i=1}^N \text{length}(l_i)\sum_{\lambda \in U_i} \text{dist}(\lambda,P_i)|\langle \phi,\phi_\lambda\rangle|^2,
$$
where $l_i$, $i=1,...,N$ are the line segments constituting the boundary of $\Omega$, $n_i$ is the unit vector orthogonal to $l_i$ directed outside $\Omega$, $P_i=\{w\in \R^2 \,|\,w\cdot n_i=0 \}$, $U_i=\{\lambda\in \Lambda \,|\,
\lambda\cdot n_i \ge 0 \}$.

Now we are ready to present the proof of Theorem \ref{BF_and_PF_invariance}.

\noindent
{\bf Theorem 1.1.}
{\it 
Let $\mu$ be the projective metaplectic representation of $SL(2,\R)$ acting on
$L^2(\R). $Let $\phi$ be a generating function of a tight Gabor frame $\{\phi_\lambda\}_{\lambda \in \Lambda}$ satisfying condition $\Phi$ with respect to $\Lambda$. Then for any $A\in SL(2,\R)$ $\Gamma =A(\Lambda)$ is a lattice, $\{(\mu(A)\phi)_\gamma\}_{\gamma \in \Gamma}$ is a tight Gabor frame with the generating function $\mu(A)\phi$ satisfying condition $\Phi$
with respect to $\Gamma$, and
\newline
(i) for any Gabor multiplier $G_{\phi,b}$ of localization type with the symbol bounded above by $1$ , $G_{\mu(A)\phi,b\circ A^{-1}}$ is a Gabor multiplier of localization type with the symbol bounded above by $1$, and
$$
PF\left(G_{\phi,b}\right)=PF\left(G_{\mu(A)\phi,b\circ A^{-1}}\right),
$$
\newline
(ii) for any $\Omega$ a $\Lambda$ domain, $A(\Omega)$ is a $\Gamma$ domain, and 
$$
BF(\phi,\Omega)=BF(\mu(A)\phi,A(\Omega)).
$$
}

\noindent
{\bf Proof.} The transformation rule $(ii)$ of Lemma \ref{transformation_properties_affine_action}, describing the effect of the conjugation by the projective metalplectic representation,  shows that
$(i)$, i.e. formula (\ref{pf_tranformation_property_gabor_multiplier})
$$
PF\left(G_{\phi,b}\right)=PF\left(G_{\mu(A)\phi,b\circ A^{-1}}\right),
$$
holds.

In the remaining part of the proof we deal with $(ii)$, i.e. with formula (\ref{transformation_property_boundary_form})
$$
BF(\phi,\Omega)=BF(\mu(A)\phi,A(\Omega)).
$$
The group $SL(2,\R)$ acts transitively on the collection of all lattices $\Gamma\subset \R^2$ satisfying condition $A_\Gamma=A_\Lambda$.
Therefore, we may assume that $A(\Lambda)=a\Z\times b\Z$, where $a,b>0$ and $ab=A_\Lambda$. It is clear that for any $A\in SL(2,\R)$ $A(\Omega)$ is an $A(\Lambda)$ domain. The proof of the $BF$ invariance rule 
(\ref{transformation_property_boundary_form}), describing the transition from $\Lambda$ to $a\Z\times b\Z$ via $A$, requires two ingredients, the 
transformation rule for the reproducing kernel 
\begin{equation}
\langle \phi, \phi_\lambda \rangle,
\label{reproducing_kernel_lambda}
\end{equation}
and for the geometric atoms
\begin{equation}
a_{s,\lambda}=\text{length}(s)\cdot \text{dist}(\lambda,l),
\label{geometric_atom_lambda}
\end{equation}
where $\lambda \in \Lambda$, $s$ is a line segment with its endpoints being 
the points of lattice $\Lambda$, constituting a portion of the boundary of $\Omega$, $l$ is the line parallel to $s$, and passing through $(0,0)$.
These two, (\ref{reproducing_kernel_lambda}) and (\ref{geometric_atom_lambda}), constitute the basic building blocks of the boundary form. Lemma \ref{transformation_properties_affine_action} $(i)$ provides the formula describing the transformation of the reproducing kernel (\ref{reproducing_kernel_lambda}) under the action of the projective metaplectic representation, it is enough to apply it. However, we need to deal with the geometric atoms (\ref{geometric_atom_lambda}), via a direct computation. 

The linear map $A \in SL(2,\R)$ providing the transition from $\Lambda$ to $a\Z \times b\Z$ satisfies condition $\text{det}(A)=1$, therefore it preserves the orientation of the boundary. It transforms a line segment $l_i$ of the boundary of $\Omega$ to the corresponding line segment $A(l_i)$ of $A(\Omega)$, and the half lattice $U_i$ of $\Lambda$ to the corresponding half lattice $A(U_i)$ of $A(\Lambda)$. The fact that $A$ might not preserve orthogonality does not matter. Lemma \ref{transformation_properties_affine_action} $(i)$ allows us to conclude that the transformation of the boundary form $BF$ is correct as far as the reproducing kernel (\ref{reproducing_kernel_lambda}) is concerned. It also shows that condition $\Phi$ transforms properly. We need to show that the geometric atoms $a_{s,\lambda}$ of (\ref{geometric_atom_lambda}) also transform properly, i.e.
that
\begin{equation}
a_{s,\lambda}=
\text{length}(As)\cdot \text{dist}(\gamma,Al),
\label{geometric_atom_azbz}
\end{equation}
where $\gamma=A\lambda$. 
The Iwasawa decomposition of 
$SL(2,\R)$ (see \cite{Lang}) allows us to represent $A$ as $PK$, where 
$P$ is upper triangular and $K$ is a rotation matrix, both in $SL(2,\R)$. Since $K$
is unitary it preserves distances and angles. It is clear that the geometric atoms
properly transform under $K$. It is therefore enough to consider the upper triangular $P$. In order to 
show that (\ref{geometric_atom_azbz}) holds,
we consider $R$, the inverse of $P$, providing a transition from $a\Z \times b\Z$
to $\Lambda$.  We verify that
\begin{equation}
\text{length}(Rs_{an,bm})\cdot \text{dist}(R(ak,bl),Rl_{an,bm})=ab|nl-km|,
\label{transformation_rules_upper_triangular_matrix}
\end{equation}
where 
$
R = \left[
	\begin{array} {cc}
		d & e \\
		0 & d^{-1}
	\end{array} \right]
$, $d>0$, $e\in \R$, $s_{an,bm}$ is the line segment connecting points $(0,0)$
and $(an,bm)$, $l_{an,bm}$ is the line containing $s_{an,bm}$, and $n,m,k,l\in \Z$. Translation invariance of the Euclidean length allows us to assume that the line segment under consideration, constituting a portion of the boundary of the 
domain of localization $\Omega$, starts at the origin $(0,0)$. 

\setlength{\unitlength}{1mm}
\begin{picture}(60,40)(-40,0)
\put(0,0){\line(2,1){60}}
\put(35,40){\circle*{1}}
\put(10,5){\circle*{1}}
\put(30,15){\circle*{1}}
\put(35,40){\line(1,-2){9}}
\thicklines
\put(10,5){\line(2,1){20}}
\put(20,5){$s_{an,bm}$}
\put(10,0){(0,0)}
\put(30,10){$(an,bm)$}
\put(40,40){$(ak,bl)$}
\put(50,19){$l_{an,bm}$}
\end{picture}
$$
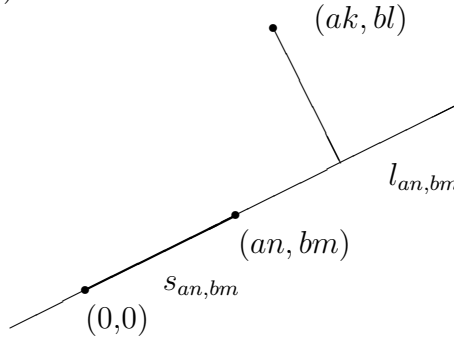
\text{{\bf Figure 1.} Displays the components of geometric atoms.}
$$
\noindent
We observe that
$R(an,bm)=(dan+ebm,bmd^{-1})$, 
$$
\text{length}^2(Rs_{an,bm})=
(dan+ebm)^2+(bmd^{-1})^2,
$$ 
the normalized vector $v$ orthogonal to $Rl_{an,bm}$ equals 
$$
v=\left( (dan+ebm)^2+(bmd^{-1})^2\right)
^{-1/2} (-bmd^{-1}, dan+ebm),
$$ 
therefore we obtain 
\begin{align}
\text{length}^2&(Rs_{an,bm})\text{dist}^2(R(ak,bl),Rl_{an,bm})=
\nonumber \\
&=\left(\langle (dak+ebl,bld^{-1}),(-bmd^{-1},dan+ebm)\rangle \right)^2
\nonumber \\
&=\left(-(dak+ebl)bm+(dan+ebm)bl \right)^2d^{-2}
\nonumber \\
&=\left( dab(nl-km)+eb^2(ml-ml)\right)^2d^{-2}=a^2b^2(nl-km)^2.
\nonumber
\end{align}
The above calculation verifies (\ref{transformation_rules_upper_triangular_matrix}). It also shows (\ref{geometric_atom_azbz}), since plugging $R=I$,
i.e. taking $d=1$, $e=0$ produces the geometric atom on $a\Z\times b\Z$.

Let us recall that Gabor multipliers with symbols of the form $b=\chi_\Omega$ are called Gabor localization operators and that they are denoted $G_{\phi,\Omega}$. The proof of Theorem \ref{BF_as_PF_limit}, the principal result of the current paper, is very lengthy and it makes use of all the auxiliary facts collected in Section 4. In order to facilitate the reading, we list the steps 
of the proof first, and then we present their proofs.

\noindent
{\bf Theorem 1.2.}
{\it 
Let $\phi$ be a generating function of a tight Gabor frame parametrized by lattice 
$\Lambda$ and satisfying condition $\Phi$, and $\Omega$ a $\Lambda$ domain 
contained in $\R^2$. Then
$$
\lim_{R \rightarrow \infty}\frac{PF(G_{\phi, R\Omega})}{R}=
\frac{1}{\text{A}_\Lambda}
BF(\phi, \Omega).
$$
}

\noindent
{\bf Steps of Proof:}
\newline
\noindent
{\bf Step 1.} We make a transition from lattice $\Lambda$ to $\Z^2$, where 
the computation is easier to handle. We construct a matrix $A\in SL(2,\R)$ and scaling parameters $a,b>0$, such that 
$\text{A}_\Lambda = ab$ and $A(\Lambda)=S_{a,b}(\Z^2)$, where $S_{a,b}$
is the scaling matrix with numbers $a,b$ on its diagonal and zeros elsewhere.
We will work with the $\Z^2$ domain $S_{a,b}^{-1}(A(\Omega))$ 
instead of the $\Lambda$ domain $\Omega$ and with the generating function
$\mu(A)\phi$ of the tight Gabor frame $\{(\mu(A)\phi)_\gamma\}_
{\gamma \in S_{a,b}(\Z^2)}$ instead of the generating function $\phi$ of the
tight Gabor frame $\{\phi_\lambda\}_{\lambda \in \Lambda}$.
Symbol $\mu$ denotes the projective metaplectic representation defined on 
$SL(2,\R)$.
We substitute the Gabor multiplier $G_{\phi, \Omega}$ by the operator
$L_{\mu(A)\phi, S_{a,b}^{-1}(A(\Omega))}$,
which is unitary equivalent to it, and defined by the formula
\begin{equation}
L_{\mu(A)\phi, S_{a,b}^{-1}(A(\Omega))}f=
\sum_{\nu \in \Z^2\cap S_{a,b}^{-1}(A(\Omega))}
\langle f,(\mu(A)\phi)_{S_{a,b}(\nu)}\rangle(\mu(A)\phi)_{S_{a,b}(\nu)}.
\label{gabor_multiplier_z_2}
\end{equation}
Once the unitary equivalence is verified it is clear that 
\begin{equation}
PF(G_{\phi, \Omega})=PF(L_{\mu(A)\phi, S_{a,b}^{-1}(A(\Omega))}).
\label{pf_g_l_equal}
\end{equation}
We substitute the boundary form $BF(\phi,\Omega)$ defined in (\ref{definition_bf}) by the boundary form $SF(\mu(A)\phi,S_{a,b}^{-1}(A(\Omega)))$ defined in terms 
of the $\Z^2$ domain $S_{a,b}^{-1}(A(\Omega)))$ and the image of the original 
generating function $\phi$ of the tight Gabor frame $\{\phi_\lambda\}_{\lambda \in \Lambda}$ under $\mu(A)$,
\begin{equation}
SF(\mu(A)\phi,S_{a,b}^{-1}(A(\Omega)))=
\sum_{i=1}^N \text{length}(l_i)\sum_{\nu \in U_i} \text{dist}(\nu,P_i)|
\langle \mu(A)\phi,(\mu(A)\phi)_{S_{a,b}(\nu)}\rangle|^2,
\label{definition_sf}
\end{equation}
where $l_i$, $i=1,...,N$ are the line segments constituting the boundary of $S_{a,b}^{-1}(A(\Omega))$, $n_i$ is the unit vector orthogonal to $l_i$ directed outside $S_{a,b}^{-1}(A(\Omega))$, $P_i=\{w\in \R^2 \,|\,w\cdot n_i=0 \}$, $U_i=\{\nu\in \Z^2 \,|\,\nu\cdot n_i \ge 0 \}$. We show that in view of condition
$\Phi$ the form (\ref{definition_sf}) is well defined and that 
\begin{equation}
SF(\mu(A)\phi,S_{a,b}^{-1}(A(\Omega)))=\frac{1}{\text{A}_\Lambda}
BF(\phi,\Omega).
\label{sf_vs_bf}
\end{equation}

\noindent
{\bf Step 2.} We use tight Gabor frames versions of Toeplitz and Hankel operators in order to express the projection functional 
$$
PF\left(L_{\mu(A)\phi, S_{a,b}^{-1}(A(R\Omega))}\right)
$$
as the square of the Hilbert-Schmidt norm of the matrix 
\begin{equation}
M_R(\nu_1,\nu_2)=\chi_{S_{a,b}^{-1}(A(R\Omega))^c}(\nu_1)\langle 
(\mu(A)\phi)_{S_{a,b}(\nu_2)}, 
(\mu(A)\phi)_{S_{a,b}(\nu_1)} \rangle \chi_{S_{a,b}^{-1}(A(R\Omega))}(\nu_2),
\label{general_kernel_form}
\end{equation}
where $\nu_1, \nu_2\in \Z^2$. As the result we obtain the equality 
\begin{equation}
\frac{PF\left(L_{\mu(A)\phi, S_{a,b}^{-1}(A(R\Omega))}\right)}{R}
=\frac{1}{R} \sum_{\nu_1,\,\nu_2 \in \Z^2} \chi_{S_{a,b}^{-1}(A(R\Omega))^c}(\nu_1)F(\nu_1-\nu_2)
\chi_{S_{a,b}^{-1}(A(R\Omega))}(\nu_2),
\label{general_kernel_form_1}\end{equation}
where 
\begin{equation}
F(\nu)=|\langle \mu(A)\phi,(\mu(A)\phi)_{S_{a,b}(\nu)} \rangle |^2.
\label{reproducing_kernel_z2}
\end{equation}
Next, we change variables and let them range over the dilated lattice $\frac{1}{R}\Z^2$. The right hand side of 
($\ref{general_kernel_form_1}$) becomes
\begin{equation}
\frac{1}{R} \sum_{\nu_1,\,\nu_2 \in \frac{1}{R}\Z^2} 
\chi_{S_{a,b}^{-1}(A(\Omega))^c}(\nu_1)F(R(\nu_1-\nu_2))
\chi_{S_{a,b}^{-1}(A(\Omega))}(\nu_2).
\label{general_kernel_form_2}\end{equation}

\noindent
{\bf Step 3.} We split the boundary of $S_{a,b}^{-1}(A(\Omega))$ into its component lattice cycles, then into their individual line segments, and then we reduce the computation of the limit to the sum of the limits over the line segments constituting the boundary. We may assume that each boundary segment connects lattice points and does not have lattice points in its interior. If 
necessary we divide it into subsegments.
Each sufficiently small neighborhood of a boundary segment gets represented as an unbounded vertical strip domain, for which expression ($\ref{general_kernel_form_2}$) becomes 
\begin{equation}
\frac{1}{R} \sum_{x_1,\,x_2\in \Z} \chi_{[0,n)}\left(\frac{x_1}{R}\right)
\chi_{[0,n)}\left(\frac{x_2}{R}\right)
\sum_{\{y_1\,|\,\frac{y_1}{R}>H\left(\frac{x_1}{R}\right)\}}
\sum_{\{y_2\,|\,\frac{y_2}{R}\le H\left(\frac{x_2}{R}\right)\}}F(x_1-x_2,y_1-y_2),
\label{vertical_strip_kernel}\end{equation}
and function $H$ representing a given segment of the boundary is a linear function with rational slope of the form $H(x)=\frac{m}{n}\,x$, with $n \ge m$, $m\ge 0$, $n>0$, and $m,n$ relatively prime. The graph of $H$ 
constitutes the top portion of the boundary
of the unbounded vertical strip domain $\{(x,y)\in \R^2 \,|\, 0\le x<n,\, y\le H(x)\}$.
Points $\left(\frac{x_1}R,\frac{y_1}R\right)$ satisfying condition 
$\frac{y_1}{R}>H\left(\frac{x_1}{R}\right)$ lie above the graph of $H$, and 
points $\left(\frac{x_2}R,\frac{y_2}R\right)$ for which 
$\frac{y_2}{R}\le H\left(\frac{x_2}{R}\right)$ 
are located  below or on the graph of $H$. Function $F$ is obtained out of
$|\langle \mu(A)\phi,(\mu(A)\phi)_{S_{a,b}(\nu)} \rangle |^2$ via a lattice transformation 
of Lemma \ref{invariance_properties} that brings the selected fragment of the boundary to the form described above. 

\noindent
{\bf Step 4.} We adjust the form of variables $x_1,x_2$ to the arithmetic form of
the slope of $H$. We represent $x_1,x_2$ as $x_1=nk_1+r_1$, $x_2=nk_2+r_2$, 
where $k_1,k_2$ are integers and $r_1,r_2 \in \{0,1,2,...,n-1\}$.  With new variables 
$k_1,r_1,k_2,r_2$ expression ($\ref{vertical_strip_kernel}$) becomes 
\begin{align}
\frac{1}{R} \sum_{r_1,\,r_2=0}^{n-1}
&\sum_{k_1,\,k_2\in \Z} 
\chi_{[0,n)}\left(\frac{nk_1+r_1}{R}\right)
\chi_{[0,n)}\left(\frac{nk_2+r_2}{R}\right)
\label{arithmetic_kernel}\\
&\sum_{\{y_1\,|\,\,y_1>mk_1+\frac{mr_1}{n}\}}
\sum_{\{y_2\,|\,\,y_2\le mk_2+\frac{mr_2}{n}\}}
F(n(k_1-k_2)+r_1-r_2,y_1-y_2).
\nonumber\end{align}
\newline
\noindent
{\bf Step 5.} Invariance of expression ($\ref{arithmetic_kernel}$) with respect 
to variables $y_1,y_2$ allows us to substitute the double summation with 
a single summation. We count the number of repetitions in the representation
$y=y_1-y_2$ and we place an appropriate factor that compensates them.
Expression ($\ref{arithmetic_kernel}$) becomes
\begin{align}
\frac{1}{R} \sum_{r_1,\,r_2=0}^{n-1}
&\sum_{k_1,\,k_2\in \Z} 
\chi_{[0,n)}\left(\frac{nk_1+r_1}{R}\right)
\chi_{[0,n)}\left(\frac{nk_2+r_2}{R}\right)
\label{invariance_kernel}\\
&\sum_{i\ge 0}
(i+1)F\left(n(k_1-k_2)+r_1-r_2,m(k_1-k_2)+\left[\frac{m}{n}r_1+1\right]-
\left[\frac{m}{n}r_2\right]+i\right),
\nonumber\end{align}
where square brackets denote the integer part of a rational number.

\noindent
{\bf Step 6.} Form ($\ref{invariance_kernel}$) is convenient for making 
the passage to the limit. It occurs that taking the limit in ($\ref{invariance_kernel}$)
with respect to $R$ is just the same as performing summation with respect to $k$. We obtain
\begin{equation}
\sum_{r_1,\,r_2=0}^{n-1}
\sum_{k \in \Z} 
\sum_{i\ge 0}
(i+1)F\left(nk+r_1-r_2,mk+\left[\frac{m}{n}r_1+1\right]-
\left[\frac{m}{n}r_2\right]+i\right).
\label{limit_kernel}\end{equation}
\newline
\noindent
{\bf Step 7.} We write expression ($\ref{limit_kernel}$) in the form
involving function $R_t(s)$
\begin{equation}
\sum_{t=0}^{n-1}
\sum_{s=0}^{n-1}
\sum_{k \in \Z} 
\sum_{i\ge 0}
(i+1)F\left(kn+t,km+R_t(s)+i\right),
\label{rts_limit_form}\end{equation}
where
\begin{equation}
R_t(s)=
\begin{cases}\left[\frac{m}{n}(t+s)+1\right]-
\left[\frac{m}{n}s\right]\text{ for }s\in\{0,1,2,...,n-t-1\}
\\
m+\left[\frac{m}{n}(s-(n-t))+1\right]-\left[\frac{m}{n}s\right]
\text{ for }s\in\{n-t,n-t+1,...,n-1\}
\end{cases}
\label{function_R_t}\end{equation}
\newline
\noindent
{\bf Step 8.} We interpret expression ($\ref{rts_limit_form}$) geometrically. 
We make use of the structural features of lattice lines and we show that in fact ($\ref{rts_limit_form}$) equals
\begin{equation}
\sqrt{n^2+m^2} \sum_{(x,\,y)\in U} \text{dist}\left((x,y),G\right)F(x,y),
\label{geometric_kernel}\end{equation}
where $G=\{(x,H(x))\,|\,x \in \R\}$ is the graph of $H$, $U$ is the part of
$\Z^2$ lying above $G$ and $\text{dist}$ is the Euclidean distance inside the 
plane $\R^2$ containing lattice $\Z^2$. Constant $\sqrt{n^2+m^2}$ represents the 
contribution of vertical strip domain $\{(x,y)\in \R^2 \,|\, 0\le x<n,\, y\le H(x)\}$
to the total length of the boundary of $S_{a,b}^{-1}(A(\Omega))$.

\noindent
{\bf Step 9.} We put together the boundary forms 
(\ref{geometric_kernel}) of vertical strip domains corresponding to all line segments of the boundary of $S_{a,b}^{-1}(A(\Omega))$ and we obtain
\begin{equation}
\lim_{R\rightarrow \infty}
\frac{PF\left( L_{\mu(A)\phi, S_{a,b}^{-1}(A(R\Omega))}\right)}{R}
=SF(\mu(A)\phi,S_{a,b}^{-1}(A(\Omega))).
\label{pf_for_l_eqls_sf}
\end{equation}
In view of (\ref{pf_g_l_equal}), and (\ref{sf_vs_bf}), formula (\ref{pf_for_l_eqls_sf}) concludes the proof. 

\noindent
{\bf Proofs of Steps:}
\newline
\noindent
{\bf Proof of Step 1.} We construct matrix $A\in SL(2,\R)$, satisfying property 
$A(\Lambda)=a\Z \times b\Z$, $a,b>0$, $ab=\text{A}_\Lambda$, by assigning values $a(0,1)$, $b(1,0)$ to the generators of  $\Lambda$ and then extending 
the assignment by linearity. The inverse of the scaling matrix $S_{a,b}$ provides a transition from $a\Z \times b\Z$
to $\Z^2$, $S_{a,b}(\Z^2)=a\Z\times b\Z$. Matrices $A$ and $S_{a,b}$ are 
the arithmetic tools needed for the transition from $\Lambda$ to $\Z^2$. The unitary operator $\mu(A)$, where $\mu$ is the projective metaplectic representation
defined on $SL(2,\R)$,
is the analytic tool responsible for the transition needed on the level of the generating functions of tight Gabor frames. 
We observe that in view of Lemma \ref{transformation_properties_affine_action}
$(ii)$ and definition (\ref{gabor_multiplier_z_2}) of $L_{\mu(A)\phi, S_{a,b}^{-1}(A(\Omega))}$
$$
\mu(A)G_{\phi,\Omega}\mu(A)^*=G_{\mu(A)\phi, A(\Omega)}
=L_{\mu(A)\phi, S_{a,b}^{-1}(A(\Omega))}.
$$
We conclude that (\ref{pf_g_l_equal}) holds.
The transition from $\Lambda$ to $\Z^2$ is justified as far the operator properties, i.e. the values of the projection functional, are concerned.

The justification of the transition from $\Lambda$ to $\Z^2$ as far as the geometric properties, i.e. the boundary form, are concerned, follows from 
Theorem \ref{BF_and_PF_invariance} $(ii)$. We obtain
$$
BF(\phi,\Omega)=BF(\mu(A)\phi,A(\Omega)),
$$
but we still need to switch to $SF(\mu(A)\phi,S_{a,b}^{-1}(A(\Omega)))$ 
defined in (\ref{definition_sf}). In order to do that, we observe that 
the boundary form scales with respect to the action of $S_{a,b}$ via the area factor $ab$, and this allows us to finish the proof of (\ref{sf_vs_bf})
$$
SF(\mu(A)\phi,S_{a,b}^{-1}(A(\Omega)))=\frac{1}{ab}BF(\phi,\Omega).
$$

\noindent
{\bf Proof of Step 2.} Since
$$L_{\mu(A)\phi,S_{a,b}^{-1}(A(R\Omega))}=G_{\mu(A)\phi,A(R\Omega)},
$$
Lemma \ref{functional_calculus} implies that 
$$
PF\left(L_{\mu(A)\phi,S_{a,b}^{-1}(A(R\Omega))}\right)=\sum_{\nu_1,\nu_2\in \Z^2}
|M_R(\nu_1,\nu_2)|^2,
$$
where $M_R(\nu_1,\nu_2)$ is defined in (\ref{general_kernel_form}). 

\noindent
{\bf Proof of Step 3.} This is the most tedious step of the proof. For the sake of notational convenience we will use symbol $\Omega$ for the current domain under consideration. It will be clear from the context what it is at a given stage of the proof. Initially $\Omega$ denotes
$S_{a,b}^{-1}(A(\Omega))$.
Our first target is to cut the kernel
$$
K_R(\nu_1,\nu_2)=\frac{1}{R}\chi_{\Omega^c}(\nu_1)F(R(\nu_1-\nu_2))
\chi_{\Omega}(\nu_2)
$$
into pieces with the help of a partition of unity of $\R^2$ obtained out of 
quadrilaterals formed around the boundary segments constituting $\partial \Omega$ and two open sets isolated from the boundary and representing the exterior and the interior of $\Omega$. 

\setlength{\unitlength}{1mm}
\begin{picture}(60,40)(-40,0)
\put(40,40){\circle*{1}}
\put(40,40){\line(2,-1){20}}
\put(60,30){\circle*{1}}
\put(60,30){\line(-1,-2){15}}
\put(45,0){\circle*{1}}
\put(45,0){\line(-5,1){15}}
\put(30,3){\circle*{1}}
\put(30,3){\line(-1,3){3}}
\put(27,12){\circle*{1}}
\put(27,12){\line(-3,-2){18}}
\put(9,0){\circle*{1}}
\put(9,0){\line(-1,2){9}}
\put(0,18){\circle*{1}}
\put(0,18){\line(1,2){11}}
\put(11,40){\circle*{1}}
\put(11,40){\line(1,0){29}}
\put(30,30){\circle*{1}}
\put(30,30){\line(1,0){15}}
\put(45,30){\circle*{1}}
\put(45,30){\line(-1,-3){5}}
\put(40,15){\circle*{1}}
\put(40,15){\line(-3,1){15}}
\put(25,20){\circle*{1}}
\put(25,20){\line(-2,-1){10}}
\put(15,15){\circle*{1}}
\put(15,15){\line(1,2){10}}
\put(25,35){\circle*{1}}
\put(25,35){\line(1,-1){5}}
\put(30,22){$\Omega^c$}
\put(38,7){$\Omega$}
\put(55,7){$\Omega^c$}
\end{picture}
$$
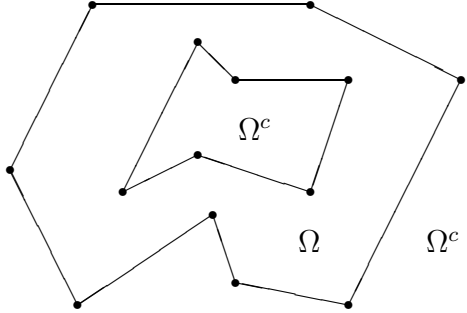
\text{{\bf Figure 2.} Explains the process of traversing the boundary $\partial \Omega$.}
$$
\noindent
We traverse each component of the
boundary of $\Omega$ according to the orientation, keeping the interior on the right and the exterior on the left, and at each initial lattice point of the boundary segment $l_i$ we place a sufficiently small rational line segment $s_i$
with its middle being the initial lattice point. We chose $s_i$ in such a way that it is transversal to both $l_i$ and the boundary segment $l_l$ preceding it,
and that both exterior and interior angles with $l_i$ and $l_l$ are smaller than $\pi$.
We assume that segments $s_i, i=1,...,N$ are so small so that they intersect the boundary $\partial \Omega$ exactly at one point. We assume that they are positioned
in such a way that
one of the endpoints of $s_i$ is inside the interior of $\Omega$ and the other in the exterior of $\Omega$.
Let $l_r$ be the segment following $l_i$. We form quadrilaterals $W_i$ out of consecutive segments $s_i,s_r$ attached at the beginning and at the end of $l_i$, and the segments joining their endpoints, both endpoints inside the interior of $\Omega$ or both in the exterior of $\Omega$.  We also assume that the sizes of the transversal segments are so small so that the segments joining the endpoints of $s_i$ and 
$s_r$ are contained either in the exterior or in the interior of $\Omega$. 
We include $s_i$ inside $W_i$, but not $s_r$. Adding $s_r$ to $W_i$ we would obtain
a closed quadrilateral, but we need to form a family of pairwise disjoint sets, and
$s_r$ is already included in $W_r$.  

\setlength{\unitlength}{1mm}
\begin{picture}(60,40)(-40,0)
\thicklines
\put(0,0){\circle*{1}}
\put(0,0){\line(2,1){20}}
\put(20,10){\circle*{1}}
\put(20,10){\line(1,1){20}}
\put(40,30){\circle*{1}}
\put(40,30){\line(2,-1){20}}
\put(60,20){\circle*{1}}
\thinlines
\put(35,40){\line(1,-2){10}}
\put(10,15){\line(2,-1){20}}
\put(30,5){\line(1,1){15}}
\put(10,15){\line(1,1){25}}
\put(22,4){$s_i$}
\put(39,35){$s_r$}
\put(39,9){$W_i$}
\put(10,1){$l_l$}
\put(31,16){$l_i$}
\put(49,27){$l_r$}
\end{picture}
$$
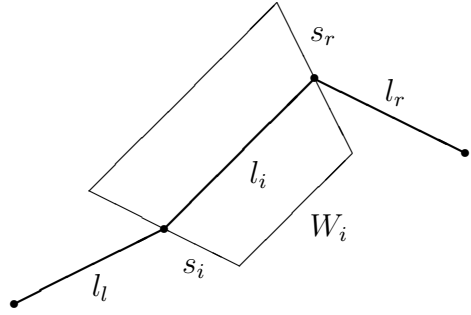
\text{{\bf Figure 3.} Illustrates the construction of quadrilaterals $W_i$.}
$$
\noindent
We define a partition of $\R^2$ out of quadrilaterals $W_i, i=1,...,N$ 
forming a neighborhood of the boundary $\partial \Omega$ and two open sets
$U_E= \Omega^c \setminus \bigcup_{i=1}^N W_i$,
$U_I= \Omega \setminus \bigcup_{i=1}^N W_i$ representing the exterior and 
the interior of $\Omega$.
Let $m_i=\chi_{W_i}, i=1,...,N$, $m_E=\chi_{U_E}$, $m_I=\chi_{U_I}$.
We define
\begin{align}
K^I_R(\nu_1,\nu_2)&=\sum_{k,l=1}^Nm_k(\nu_1)K_R(\nu_1,\nu_2)m_l(\nu_2),\nonumber \\
K^{II}_R(\nu_1,\nu_2)&=\sum_{l=1}^Nm_E(\nu_1)K_R(\nu_1,\nu_2)m_l(\nu_2),\nonumber \\
K^{III}_R(\nu_1,\nu_2)&=\sum_{k=1}^Nm_k(\nu_1)K_R(\nu_1,\nu_2)m_I(\nu_2),\nonumber \\
K^{IV}_R(\nu_1,\nu_2)&=m_E(\nu_1)K_R(\nu_1,\nu_2)m_I(\nu_2).
\nonumber \end{align}
We observe that
\begin{equation}
K_R=\sum_{J=I}^{IV}K^J_R.
\label{decomposition_of_K}\end{equation}
We do not need terms with factors $m_I(\nu_1)$ and  $m_E(\nu_2)$ 
in the above decomposition, since
$U_I\,\cap \,\Omega^c=\emptyset$ and 
$U_E\,\cap \,\Omega=\emptyset$. We know that
$\text{dist}(U_I, \Omega^c)>0$, 
$\text{dist}(U_E, \Omega)>0$, therefore Lemma \ref{separated_supports} allows us to conclude
\begin{equation}
\lim_{R\rightarrow \infty}\sum_{\nu_1,\nu_2 \in \frac{1}{R}\Lambda}
K^J_R(\nu_1, \nu_2)=0,
\label{eliminate_I_E} \end{equation}
for $J=II,III,IV$. 
Since $\text{dist}(W_k, W_l)>0$
for $k\ne l$ such that $l_k,l_l$ are not neighboring line segments, again
by  Lemma \ref{separated_supports} we obtain
\begin{equation}
\lim_{R\rightarrow \infty}\sum_{\nu_1,\nu_2 \in \frac{1}{R}\Lambda}
m_k(\nu_1)K_{R}(\nu_1, \nu_2)m_l(\nu_2)=0.
\label{eliminate_non_neighbors} \end{equation}
If $k \ne l$, but $l_k, l_l$ are neighboring line segments, then Lemma
 \ref{supports_in_cones} applies and we conclude again 
\begin{equation}
\lim_{R\rightarrow \infty}\sum_{\nu_1,\nu_2 \in \frac{1}{R}\Lambda}
m_k(\nu_1)K_{R}(\nu_1, \nu_2)m_l(\nu_2)=0.
\label{eliminate_neighbors} \end{equation}

\setlength{\unitlength}{1mm}
\begin{picture}(60,40)(-40,0)
\thicklines
\put(10,10){\circle*{1}}
\put(10,10){\line(1,1){20}}
\put(30,30){\circle*{1}}
\put(30,30){\line(1,0){25}}
\put(55,30){\circle*{1}}
\thinlines
\put(25,40){\line(1,-2){10}}
\put(0,15){\line(2,-1){20}}
\put(20,5){\line(1,1){15}}
\put(0,15){\line(1,1){25}}
\put(60,40){\line(-1,-2){10}}
\put(25,40){\line(1,0){35}}
\put(50,20){\line(-1,0){15}}
\put(40,15){$W_l$}
\put(28,10){$W_k$}
\put(21,16){$l_k$}
\put(40,32){$l_l$}
\put(54,23){$W_l\cap \Omega$}
\put(40,23){$\Omega$}
\put(0,31){$W_k\cap \Omega^c$}
\put(14,22){$\Omega^c$}
\end{picture}
$$
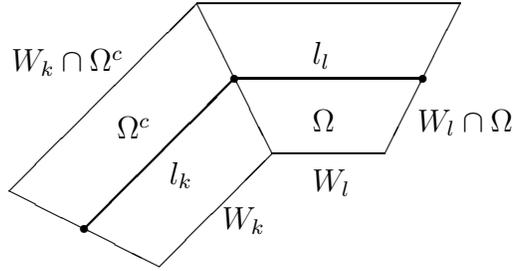
\text{{\bf Figure 4.} Explains the usage of Lemma \ref{supports_in_cones}
for neighboring segments $l_k, l_l$.}
$$
\noindent
Formulae  (\ref{decomposition_of_K})-(\ref{eliminate_neighbors})
allow us to conclude that only terms of the form
$$
m_k(\nu_1)K_{R}(\nu_1, \nu_2)
m_k(\nu_2)
$$
contribute to the limit.

Let us recall that we assumed that each boundary segment $l_k$ connects lattice points and does not have lattice points in its interior. We may apply Lemma \ref{invariance_properties} in order to bring 
each segment  $l_k$ to the form $(x,H(x))$, 
with $0\le x< n$, $H(x)=\frac{m}{n}x$, $n>0$, $m\ge 0$, $n\ge m$, $n,m$ relatively prime, with the image of $\Omega$ placed below the graph of $H$, and the 
image of  
$\Omega^c$ placed above the graph of $H$. Under the transformation  of Lemma \ref{invariance_properties} quadrilateral $W_k$ becomes $\tilde{W}_k$. With the help of Lemma
\ref{supports_in_cones} we bring $\tilde{W}_k$ to a form of a 
bounded vertical strip domain with the top portion of the boundary
represented by a line $t_k$ and the bottom portion by a line $b_k$. Lemma \ref{separated_supports} allows us to substitute the bounded vertical strip
domain we have just obtained by a rectangle with the top boundary
represented by a horizontal line $y=r$, with $r\in \Z$, placed above $t_k$,
and the bottom boundary represented by a horizontal line $y=s$, with $s\in \Z$,
placed below $b_k$. 

\setlength{\unitlength}{1mm}
\begin{picture}(60,80)(-40,0)
\thicklines
\put(25,30){\circle*{1}}
\put(25,30){\line(1,1){20}}
\put(45,50){\circle*{1}}
\thinlines
\put(40,60){\line(1,-2){10}}
\put(15,35){\line(2,-1){20}}
\put(25,15){\line(1,1){25}}
\put(15,35){\line(1,1){30}}
\put(8,34){$\tilde{W}_k$}
\put(37,35){$\Omega$}
\put(29,42){$\Omega^c$}
\put(25,0){\line(0,1){80}}
\put(45,0){\line(0,1){80}}
\put(25,5){\line(1,0){20}}
\put(25,75){\line(1,0){20}}
\put(22,4){$s$}
\put(22,74){$t$}
\put(40,26){$b_k$}
\put(27,53){$t_k$}
\end{picture}
$$
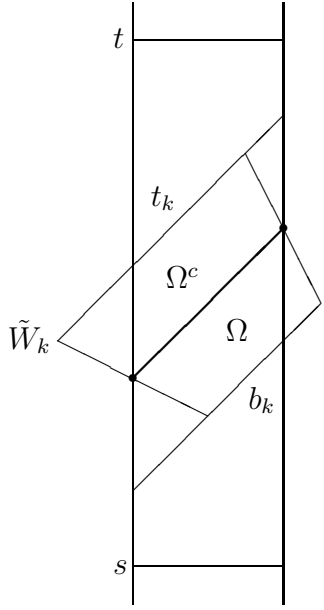
\text{{\bf Figure 5.} Illustrates the usage of quadrilateral $\tilde{W}_k$.}
$$
\noindent
In the last step we switch to the unbounded strip
domain 
$$
\Omega=\{(x,y)\,|\,0\le x <n, y\le H(x) \},
$$
with its complement
$$
\Omega^c=\{(x,y)\,|\,0\le x <n, y> H(x) \}.
$$
We have the splitting
\begin{align}
\chi_{\Omega^c}(\nu_1)F(R(\nu_1-\nu_2))\chi_{\Omega}(\nu_2)&=
\nonumber \\
\chi_{\{ (x,y)\,|\,0\le x<n, y\le r\}}(\nu_1)&\chi_{\Omega^c}(\nu_1)F(R(\nu_1-\nu_2))\chi_{\Omega}(\nu_2)
\chi_{\{ (x,y)\,|\,0\le x<n, y > s\}}(\nu_2)
\label{<r>s}\\
+\chi_{\{ (x,y)\,|\,0\le x<n, y > r\}}(\nu_1)&\chi_{\Omega^c}(\nu_1)F(R(\nu_1-\nu_2))\chi_{\Omega}(\nu_2)
\chi_{\{ (x,y)\,|\,0\le x<n, y > s\}}(\nu_2)
\label{>r>s}\\
+\chi_{\{ (x,y)\,|\,0\le x<n, y\le r\}}(\nu_1)&\chi_{\Omega^c}(\nu_1)F(R(\nu_1-\nu_2))\chi_{\Omega}(\nu_2)
\chi_{\{ (x,y)\,|\,0\le x<n, y \le s\}}(\nu_2)
\label{<r<s}\\
+\chi_{\{ (x,y)\,|\,0\le x<n, y> r\}}(\nu_1)&\chi_{\Omega^c}(\nu_1)F(R(\nu_1-\nu_2))\chi_{\Omega}(\nu_2)
\chi_{\{ (x,y)\,|\,0\le x<n, y \le s\}}(\nu_2)
\label{>r<s}
\end{align}
We need to show that the normalized sums coming out of terms (\ref{>r>s}),
(\ref{<r<s}), (\ref{>r<s}) have zero limits. Lemma \ref{separated_supports}
applies to both (\ref{>r>s}), (\ref{<r<s}), but (\ref{>r<s}) needs to be treated
separately. After the change of variables the normalized sum of (\ref{>r<s}) becomes
$$
\frac{1}{R}\sum_{0\le x_1,x_2 <Rn}\sum_{\substack{y_1>rR \\
y_2 \le sR}}F(x_1-x_2, y_1-y_2).
$$
In the proof that its limit is zero as $R\rightarrow \infty$ we may assume that $R$ is
an integer. Let $u=\text{min}\{r,-s\}$ and let $*$ denote the convolution 
on $\Z$. We have 
\begin{align}
\frac{1}{R}\sum_{0\le x_1,x_2 <Rn}\sum_{\substack{y_1>rR \\
y_2 \le sR}}&F(x_1-x_2, y_1-y_2) 
\le \frac{1}{R}\sum_{0\le x_1,x_2 <Rn}\sum_{y_1, y_2 \ge uR}F(x_1-x_2, y_1+y_2)
\nonumber \\
&\le \frac{1}{R}\sum_{0\le x_1,x_2 <Rn}\sum_{y \ge 2uR}F(x_1-x_2, y)
(y-2uR+1) 
\nonumber \\
&\le n \sum_{y\ge 2uR}\left\langle F(\cdot, y)(y+1)*\frac{\chi_{[0,Rn)}}{(Rn)^{1/2}}
,\frac{\chi_{[0,Rn)}}{(Rn)^{1/2}} \right\rangle
\nonumber \\
&\le n \sum_{y\ge 2uR}\sum_{x\in \Z} F(x,y)(y+1) \ra_{R \rightarrow \infty} 0,
\nonumber
\end{align}
because $2uR \rightarrow \infty$ as $R\rightarrow \infty$, since $u>0$, and
$$
\sum_{y\in \Z}\sum_{x\in \Z}F(x,y)(|y|+1)<\infty.
$$

\noindent
{\bf Proof of Step 4.} With the new representation of $x_1,x_2$ the summation done in (\ref{vertical_strip_kernel}) becomes (\ref{arithmetic_kernel}). It is enough to perform substitutions and then simplify the resulting expressions.

\noindent
{\bf Proof of Step 5.} We have 
\begin{align}
y_1>m k_1+\frac{m r_1}{n}&\text{, i.e. }y_1\ge m k_1+\left[\frac{m r_1}{n}+1\right],\nonumber \\
y_2\le m k_2+\frac{m r_2}{n}&\text{, i.e. }y_2\le m k_2+\left[\frac{m r_2}{n}\right],
\nonumber\end{align}
therefore 
$$
y_1-y_2= m(k_1-k_2)+\left[\frac{mr_1}{n}+1\right] - \left[\frac{mr_2}{n}\right]+i,
$$
for some $i=0,1,2,...$. A given value $i$ occurs for exactly  $i+1$ pairs
$y_1,y_2$.

\noindent 
{\bf Proof of Step 6.} We move summations over $k_1,k_2$ in (\ref{invariance_kernel}) inside and summations over $r_1,r_2,i$ outside. Our target is to identify convolution kernels defined in terms of variables $k_1,k_2$ and then, keeping variables $r_1,r_2,i$ fixed, take the limit with respect to $R$. For a rational number $q$ by $\lceil q \rceil$ we denote the smallest integer larger or equal to $q$. We know that 
$$
0\le \frac{n k_i+r_i}{R} < n, \, i=1,2.
$$
We observe that $k_i$ ranges over the integer interval
$$
\left[ 0,\lceil R-\frac{r_i}{n}\rceil -1 \right], \,i=1,2. 
$$
We define convolution kernels $K_{r_1,r_2,i}$ as  
$$
K_{r_1,r_2,i}(k)=F\left(n k+r_1-r_2,m k+\left[\frac{m}{n}r_1+1\right]-
\left[\frac{m}{n}r_2\right]+i\right).
$$
Lemma \ref{approximate_identity} guaranties that as $R\rightarrow \infty$
the expression
$$
\frac{1}{\lceil R-\frac{r_i}{n}\rceil} \sum_{k_1,\,k_2\in \Z} 
K_{r_1,r_2,i}(k_1-k_2)
\chi_{\left[0,\lceil R-\frac{r_i}{n}\rceil -1\right]}\left(k_1\right)
\chi_{\left[0,\lceil R-\frac{r_i}{n}\rceil -1\right]}\left(k_2\right)
$$
tends to $\sum_k K_{r_1,r_2,i}(k)$. We proved that (\ref{limit_kernel}) is the limit of (\ref{invariance_kernel}).

\noindent
{\bf Proof of Step 7.} Each number $kn+t$, $t\in\{1,2,...,n-1\}$
has two possible representations. The first one as $kn+r_1-r_2$
with $r_1-r_2=t$ and the second one as $(k+1)n+r_1-r_2$
with $r_1-r_2=t-n$. The first case occurs for $r_1=t+s$, $r_2=s$,
$s\in \{0,1,...n-t-1\}$ and the second case for $r_1=s-(n-t)$, $r_2=s$,
$s\in \{n-t,n-t+1,...,n-1\}$. For $t=0$ we have only one representation of 
$kn+t$ as $kn+r_1-r_2$ with $r_1=r_2=s$, $s\in \{0,1,...n-1\}$.
Function $R_t(s)$ incorporates all these relations and allows us to subsitue summations with respect to $r_1,r_2$ by summations with $t,s$, therefore
(\ref{limit_kernel}) becomes (\ref{rts_limit_form}).

\noindent
{\bf Proof of Step 8.} We fix summation variables $k,t$ of expression
(\ref{rts_limit_form}) and we consider the effect of the summation done
with respect to $s,i$. By Lemma \ref{lattice_slopes} we know that 
function $R_t(s)$ takes two values $\left[ \frac{m}{n}t+1\right]$ and
$\left[ \frac{m}{n}t+1\right]+1$, therefore we conclude that points
$(kn+t,km+R_t(s)+i)$, $s=0,1,...,n-1$, $i=0,1,2,...$, represent a discrete 
vertical half line consisting of lattice points of $\Lambda$ starting
directly above the graph $G$. Lemma  \ref{lattice_slopes} also 
tells us that value $\left[ \frac{m}{n}t+1\right]$ is taken 
$n\left(\left[ \frac{m}{n}t+1\right]-\frac{m}{n}t\right)$ times and
value  $\left[ \frac{m}{n}t+1\right]+1$ the remaining 
$n-n\left(\left[ \frac{m}{n}t+1\right]-\frac{m}{n}t\right)$ times.
The first point of the discrete half line, located  right above the graph $G$, corresponds to parameter values $i=0$, $R_t(s)=\left[ \frac{m}{n}t+1\right]$,
and it is repeated inside formula (\ref{rts_limit_form}) $n\left(\left[ \frac{m}{n}t+1\right]-\frac{m}{n}t\right)$
times. The $j$th point of it, $j\ge 2$, with counting done upwards, corresponds to values $i=j-1$, 
$R_t(s)=\left[ \frac{m}{n}t+1\right]$ and $i=j-2$, 
$R_t(s)=\left[ \frac{m}{n}t+1\right]+1$, and in (\ref{rts_limit_form}) it comes with multiplicity
$n\left( j-1+\left[ \frac{m}{n}t+1\right]-\frac{m}{n}t\right)$.
Let $d_j,v_j$ represent the Euclidean and the vertical distances
from the $j$th point, $j\ge 1$ of the discrete half line to the graph $G$.
Similarity relation of the triangle representing distances $d_i,v_i$ with the triangle with vertices $(0,0),(n,0),(n,m)$ gives
$$
\frac{d_j}{v_j}=\frac{n}{\sqrt{n^2+m^2}}.
$$
Since the multiplicity inside (\ref{rts_limit_form}) of the $j$th point of the discrete half line  equals $nv_j$, and variables $k,t$ parameterize 
all discrete vertical half lines of $U$ we conclude that (\ref{rts_limit_form}) 
may be expressed as (\ref{geometric_kernel}).

\noindent
{\bf Proof of Step 9.} It is enough to combine together the outcomes of 
steps 1-8, i.e.  all of the intermediate stages of the reduction process.

The proof of Corollary \ref{no_best_lattice} follows directly from Theorems \ref{BF_and_PF_invariance}, \ref{BF_as_PF_limit}.

\noindent
{\bf Corollary 1.3.} {\it
For any lattice $\Lambda \subset \R^2$ satisfying condition $\text{A}_\Lambda <1$, any generating function $\phi$ of a tight Gabor frame $\{ \phi_\lambda \}_{\lambda\in \Lambda}$, any $\Lambda$ lattice domain $\Omega$, and any $a,b>0$ satisfying $ab=\text{A}_\Lambda$, there are a generating function $\breve{\phi}$ of a tight Gabor frame $\{ \breve{\phi}_\lambda \}_{\lambda\in a\Z \times b\Z}$ and a $a\Z \times b\Z$ lattice domain $\breve{\Omega}$, satisfying $||\phi||_{L^2(\R)}=||\breve{\phi}||_{L^2(\R)}$, $\text{Area}( {\Omega})=\text{Area}( \breve{\Omega})$, $PF(G_{\phi,R\Omega})=
PF(G_{\breve{\phi},R\breve{\Omega}})$, for all $R>0$, and also 
$BF(\phi, \Omega)=BF(\breve{\phi}, \breve{\Omega})$.
The rates of convergence of $\frac{PF(G_{\phi,R\Omega})}{R}$ to 
$\frac{1}{A_\Lambda}BF(\phi,\Omega)$ and $\frac{PF(G_{\breve{\phi},R\breve{\Omega}})}{R}$ to 
$\frac{1}{ab}BF(\breve{\phi},\breve{\Omega})$ are the same.
}

\noindent
{\bf Proof.} It is enough to take $A\in SL(2,\R)$ transferring $\Lambda$ onto
$a\Z\times b\Z$, $\breve{\phi}=\mu(A)\phi$, $\breve{\Omega}=A(\Omega)$, and apply Theorem \ref{BF_and_PF_invariance}. 
The existence and the form of the limits follows from Theorem \ref{BF_as_PF_limit}.

\section{Auxiliary Facts and their Proofs}

\noindent
{\bf Transformation properties with respect to the metaplectic representation.} A comprehensive presentation of the metaplectic representation from the point of view of phase space analysis is contained in Folland's book \cite{Fol}. The book by Lang \cite{Lang} is an extensive reference on $SL(2,\R)$.
The Heisenberg group $\mathbb{H}^1$ is the group obtained 
by defining on $\R^{3}$ the product
$$
(z,t)\cdot (w,s)=\left(z+w,t+s-\frac{1}{2}\omega(z,w)\right),
$$
where $z,w\in \R^{2}$, $t,s\in \R$ and $\omega$ is the symplectic form defined on 
$\R^{2}$, i.e. 
$$
\omega (z,w)= z^t Jw,\,
J = \left[
	\begin{array} {cc}
		0 & 1 \\
		-1 & 0
	\end{array} \right].
$$ 
The Schr\"odinger representation of the group $\mathbb{H}^1$,
acting on $L^2(\R)$, is then defined by
$$
\rho(x,\xi,t)f(y)=e^{2\pi i t}
e^{2\pi i \xi y}f(y-x).
$$
We write $z=(x,\xi)$ when we separate the position component $x$
from the momentum component $\xi$ of a point $z$ of the phase space $\R^{2}$.
The group $SL(2,\R)$, consisting of $2\times 2$ matrices with real entries
and determinant $1$,  
acts on $\mathbb{H}^1$ via automorphisms that leave the center
$\{(0,t)\,|\,t\in \R\}$ of  $\mathbb{H}^1$ pointwise fixed,
$$
A\cdot (z,t)=(Az,t).
$$
For any fixed $A\in SL(2,\R)$ there is a unitary representation of 
$\mathbb{H}^1$, acting on $L^2(\R)$, defined as the composition
$$
\rho_A(z,t) = \rho(A\cdot (z,t)),
$$
with its restriction to the center of $\mathbb{H}^1$ being a multiple of the identity. By the Stone-von Neumann theorem $\rho_A$ is unitary equivalent to $\rho$,
i.e. there is an intertwining unitary operator $\mu(A)$ acting on $L^2(\R)$
such that for all $(z,t)\in \mathbb{H}^1$
$$
\rho_A(z,t)=\mu(A)\circ \rho(z,t)\circ \mu(A)^{-1}.
$$ 
By Schur's lemma, $\mu$
is determined up to a phase factor $e^{is}, s\in \R$. It turns out 
that the phase ambiguity is really a sign, so that $\mu$ lifts to a representation
of the double cover of the group $SL(2,\R)$. The constructed 
representation of the double cover of $SL(2,\R)$ is called the metaplectic
representation.

The representations $\rho$ and $\mu$ can be combined and give rise to the
extended metaplectic representation  $\mu_e$, the composition of 
$\rho$, defined on $\mathbb{H}^1$, with $\mu$, defined on the double 
cover of $SL(2,\R)$. From the point of view of the interpretation as a phase space action, the phase factors do not matter, therefore we remove them and treat 
$\mu$ as a projective representation of $SL(2,\R)$, and $\mu_e$ as a projective representation of the semidirect product  $\R^{2}\rtimes SL(2,\R)$
with the group law
\begin{equation}
\left(z,A\right)\cdot \left( w,B \right)=
\left( z+ Aw,AB \right).
\label{affine_phase_space_transformations_group_law}
\end{equation}
The extended metaplectic representation provides all affine transformations of the phase space $\R^{2}$. For $\left(z,A\right) \in \R^{2}\rtimes SL(2,\R)$
the unitary operator $\mu_e\left(z,A\right)$, defined up to a phase factor,
expresses the analytic action on $L^2(\R)$. The affine geometric action 
on $\R^2$ is expressed by the law
\begin{equation}
\left( z,A\right)w=Aw+z,\text{ where }w\in \R^2.
\label{affine_phase_space_transformations_group_action}
\end{equation}
The extended metaplectic representation is a convenient setup for performing computations involving compositions of the Schr\"odinger and
the metaplectic representations.

Both Gabor multipliers (\ref{gabor_multiplier}) and lattice boundary forms 
(\ref{definition_bf}) have natural transformation properties with respect to the projective metataplectic representation. These properties are important ingredient of our proofs. We will deduce them out of the  
the fundamental lemma formulated below.
\begin{lem}\label{transformation_properties_affine_action}
Let $\Lambda \subset \R^2$ be a lattice, $\{\phi_\lambda\}_{\lambda \in \Lambda}$ a tight Gabor frame with the generating function $\phi$, 
$A\in SL(2,\R)$, and $b\in l^\infty(\Lambda)$. Then
\newline
(i) $\Gamma=A\Lambda \subset \R^2$ is a lattice, $\{(\mu(A)\phi)_\gamma\}_{\gamma \in \Gamma}$ a tight Gabor frame with the generating function $\mu(A)\phi$, and
\begin{equation}
\langle \phi,\phi_\lambda \rangle = 
\langle \mu(A)\phi, (\mu(A)\phi)_\gamma \rangle
\text{ for }\gamma = A\lambda,
\label{tranformation_property_tight_frame}
\end{equation} 
\newline
(ii) $b\circ A^{-1} \in l^\infty (\Gamma)$, and 
\begin{equation}
\mu(A)G_{\phi,b}\mu(A)^*=G_{\mu(A)\phi, b\circ A^{-1}}.
\label{transformation_property_gabor_multiplier}
\end{equation}
\end{lem}

\noindent
{\bf Proof.} Clearly $\Gamma=A\Lambda$ is a lattice, since $A$ is linear
and invertible. Group law (\ref{affine_phase_space_transformations_group_law}) and affine action rule (\ref{affine_phase_space_transformations_group_action})
allow us to identify the phase space $\R^2$ with the subgroup
$\{((q,p), I)\,|\,q,p\in \R\}$ of $\R^{2}\rtimes SL(2,\R)$. Since
$(0,A)(\lambda,I)(0,A^{-1})=(A\lambda,I)$, by substituting $\gamma=A\lambda$ we obtain
\begin{align}
\mu(A)&G_{\phi,b}\mu(A)^*f =\sum_{\lambda \in \Lambda}  
b(\lambda)\langle f,\mu(A)(\phi_\lambda)\rangle \mu(A)(\phi_\lambda) 
\nonumber \\
&=  \sum_{\lambda \in \Lambda} 
b(\lambda)\langle f,\mu_e((0,A)(\lambda,I))\phi \rangle 
\mu_e((0,A)(\lambda,I))\phi\
\nonumber \\
&= \sum_{\lambda \in \Lambda}  
b(\lambda)\langle f,\mu_e((0,A)(\lambda,I)(0,A^{-1}))\mu(A)\phi \rangle \mu_e((0,A)(\lambda,I)(0,A^{-1}))\mu(A)\phi
\nonumber \\
&= \sum_{\gamma \in \Gamma} 
b(A^{-1}\gamma)\langle f,\mu_e(\gamma,I)\mu(A)\phi \rangle \mu_e(\gamma,I)\mu(A)\phi
\nonumber \\
&= \sum_{\gamma \in \Gamma} 
b(A^{-1}\gamma)\langle f,(\mu(A)\phi)_\gamma \rangle (\mu(A)\phi)_\gamma
\nonumber \\
&=G_{\mu(A)\phi, b\circ A^{-1}}f.
\nonumber \end{align}
The above calculation shows that $\{(\mu(A)\phi)_\gamma\}_{\gamma \in\Gamma}$ is a tight Gabor frame, it is enough to take the constant function equal 
to $1$ for $b$. It also verifies formula (\ref{transformation_property_gabor_multiplier}). 
Formula (\ref{tranformation_property_tight_frame}) follows by a similar calculation,
$$
\langle \mu(A)\phi, (\mu(A)\phi)_\gamma \rangle=
\langle \phi, \mu_e((0,A^{-1})(\gamma,I)(0,A))\phi\rangle=
\langle \phi,\mu_e(A^{-1}\gamma,I)\phi\rangle=
\langle \phi,\phi_\lambda\rangle.
$$

\noindent
{\bf Symbolic calculus of Gabor multipliers.} Let us assume that  $\Lambda\subset \R^2$ is a lattice, and $\{\phi_\lambda\}_{\lambda \in \Lambda}$ a tight Gabor frame defined on it. Let us define the mapping 
$W_\phi: L^2(\R)\rightarrow l^2(\Lambda)$ by the formula
$$
W_\phi f(\lambda)=\langle f,\phi_\lambda \rangle.
$$
Tight frame properties of $\{\phi_\lambda\}_{\lambda \in \Lambda}$ imply that $W_\phi$ is an isometry, and that the operator
$P_\phi: l^2(\Lambda )\rightarrow  l^2(\Lambda )$ defined as
$$
P_\phi h(\lambda)=\sum_{\mu \in \Lambda}\langle \phi_\mu,\phi_\lambda\rangle h(\mu)
$$ 
is the orthogonal
projection onto the range of $W_\phi$. Gabor multiplier $G_{\phi,b}$ is parametrized by the generating function $\phi$ of a tight Gabor frame 
$\{\phi_\lambda\}_{\lambda \in \Lambda}$
and a symbol $b\in l^\infty (\Lambda)$.  Let us recall that it is a bounded 
operator acting on $L^2(\R)$ defined as
$$
G_{\phi,b}(f)=\sum_{\lambda \in \Lambda}
b(\lambda)\langle f,
 \phi_\lambda\rangle\phi_\lambda.
$$
It is convenient to describe Gabor multiplier $G_{\phi,b}$ in 
terms of the Toeplitz operator $\mathcal{G}_{\phi,b}=P_\phi M_bP_\phi$ acting on $l^2(\Lambda)$, $M_b$ denotes the operator of multiplication by $b$. The isometry $W_\phi$ allows us to identify $G_{\phi,b}$ with the upper left corner of the matrix representation of the operator ${\mathcal{G}}_{\phi,b}$
with respect to the orthogonal decomposition $l^2(\Lambda)=W_\phi(L^2(\R))
\oplus W_\phi(L^2(\R))^\perp$. Basic properties of the symbolic calculus 
of Gabor multipliers can be deduced out of the properties of the symbolic calculus of Toeplitz operators. Hankel operator $\mathcal{H}_{\phi,b}=(I-P_\phi )M_bP_\phi$ acting on $l^2(\Lambda)$ measures to what degree the mapping 
$b\rightarrow \mathcal{G}_{\phi,b}$
fails to be a homomorphism. The algebraic formula
$$
\mathcal{G}_{\phi,b_1b_2}-\mathcal{G}_{\phi,b_1}\mathcal{G}_{\phi,b_2}=
\mathcal{H}_{\phi,\overline{b_1}}^*\mathcal{H}_{\phi,b_2}
$$
expresses this relationship quantitatively and it is the main conceptual ingredient 
of the argument that allows us to write down the projection functional 
$$
PF(G_{\phi,\Omega})=
PF(\mathcal{G}_{\phi,\Omega})=
\text{tr}\left( \mathcal{G}_{\phi,\Omega}\left( 
I-\mathcal{G}_{\phi,\Omega}\right) \right)
$$ 
as the square of the Hilbert-Schmidt norm the matrix
$$
M(\nu_1,\nu_2)
=\chi_{\Omega^c}(\nu_1)\langle \phi_{\nu_2}, \phi_{\nu_1} \rangle \chi_{\Omega}(\nu_2),
$$ 
where $\nu_1,\nu_2\in \Lambda$.
\begin{lem}\label{functional_calculus} Let $\{\phi_\lambda\}
_{\lambda \in \Lambda}$ be a tight Gabor frame and $\Omega \subset 
\Lambda$ a finite set. Then
$$
PF(G_{\phi,\Omega})=\sum_{\nu_1,\nu_2\in \Lambda}
\chi_{\Omega^c}(\nu_1)|\langle \phi_{\nu_2}, \phi_{\nu_1}\rangle |^2 \chi_{\Omega}(\nu_2).
$$
\end{lem}

\noindent
{\bf Proof.} The non-zero eigenvalues of the localization operator $G_{\phi,\Omega}$ coincide with the non-zero eigenvalues of the Toeplitz operator 
$\mathcal{G}_{\phi,\Omega}=P_\phi M_{\chi_\Omega} P_\phi$ and the non-zero 
eigenvalues of the operator $\mathcal{R}_{\phi,\Omega}=
M_{\chi_\Omega}P_\phi M_{\chi_\Omega}$.
The first fact follows from the identification of $G_{\phi,\Omega}$ with the 
upper left corner of the matrix representation of $\mathcal{G}_{\phi,\Omega}$ 
with respect to the decomposition of $l^2(\Lambda)$ into the range 
of $W_\phi$ and its orthogonal complement. The second fact follows since
for a compact operator $T$ the non-zero eigenvalues of $TT^*$ and 
$T^*T$ are the same. We observe that 
$$
PF(G_{\phi,\Omega})=PF(\mathcal{G}_{\phi,\Omega})=
\text{tr}\left( \mathcal{G}_{\phi,\Omega}\left(I-
\mathcal{G}_{\phi,\Omega}\right)\right)=
\text{tr}\left( \mathcal{R}_{\phi,\Omega}\left(I-
\mathcal{R}_{\phi,\Omega}\right)\right)
$$ 
and that 
$$
M_{\chi_\Omega}P_\phi M_{\chi_\Omega}(I-M_{\chi_\Omega}P_\phi M_{\chi_\Omega})=M_{\chi_\Omega}P_\phi (I-M_{\chi_\Omega})P_\phi M_{\chi_\Omega}=M_{\chi_\Omega}P_\phi M_{\chi_{\Omega^c}} P_\phi M_{\chi_\Omega}.
$$
Therefore we obtain
$$
PF(G_{\phi,\Omega})= \text{tr}\left(M_{\chi_\Omega}P_\phi M_{\chi_{\Omega^c}} P_\phi M_{\chi_\Omega}\right)=\sum_{\nu_1,\nu_2\in \Lambda}
\chi_{\Omega^c}(\nu_1)|\langle \phi_{\nu_2}, \phi_{\nu_1}\rangle |^2 \chi_{\Omega}(\nu_2),
$$
and this finishes the proof.

\noindent
{\bf Lattice slopes of rational lines.} For an integer $t$ we define the 
{\it $t$-slope} of the line $y=ax+b$ at an integer argument $s$ as 
\begin{equation}
S_t(s)=U(a(t+s)+b)-l(as+b),
\label{t-slope}\end{equation}
where $U(x)$ is the smallest integer larger than $x$, and $l(x)$ is the largest integer 
smaller or equal to $x$. The $t$-slope at $s$ is simply the smallest
difference between integer values above the graph at $t+s$ and below or on the graph at $s$. 
We interpret it as the lattice slope corresponding to making $t$-steps to the right of $s$.
We do not normalize the $t$-slope, i.e. we do not divide it by the number of steps. 
As we have already seen function $R_t$ defined in (\ref{function_R_t}) is the 
principal analytic component of the boundary form. It occurs that it may be interpreted as 
the $t$-slope of the rational line $y=\frac{m}{n}x$. Indeed, let $S_t(s)$ be the $t$-slope 
of the line $y=\frac{m}{n}x$ at $s$. Direct inspection shows, that values $R_t(s)$
and $S_t(s)$ coincide for $s\in\{0,1,2,...,n-1\}$. Function $S_t$ is periodic with 
period $n$ and we may regard it as defined on the cyclic group $\Z_n$.
\par
Our primary geometric concern are the values of lattice slopes and the frequencies 
with which they occur. Observe that $R_0(s)=1$ for all $s\in\{0,1,2,...,n-1\}$. 
If $n=1$, then the only possible value of $t$ is $0$ and again $1$ is the only value of 
$R_t$. If however $n\ge 2$ and $t\ne 0$, then $R_t$ takes precisely two values. The next 
lemma describes those values and the frequencies with which they occur. 
\begin{lem}\label{lattice_slopes} Let $m,n$ be a pair of relatively prime integers. 
Assume that $n\ge 2$. Let $t\in\{0,1,2,...,n-1\}$ be a fixed number. Function
$R_t(s)$ defined for arguments $s\in\{0,1,2,...,n-1\}$ by formula (\ref{function_R_t})
takes two distinct values $\left[\frac{m}{n}t+1\right]$ and $\left[\frac{m}{n}t+1\right]+1$. 
Value $\left[\frac{m}{n}t+1\right]$ is taken 
$n\left(\left[\frac{m}{n}t+1\right]-\frac{m}{n}t\right)$
times, and value $\left[\frac{m}{n}t+1\right]+1$ is taken the remaining 
$n-n\left(\left[\frac{m}{n}t+1\right]-\frac{m}{n}t\right)$ times.
\end{lem}

\noindent
{\bf Proof.} We know that $R_t$ and $S_t$ are equal. It is therefore enough 
to prove Lemma \ref{lattice_slopes} with $S_t$ instead of $R_t$.
Let $s\in\{0,1,2,...,n-1\}$. Line $y=\frac{m}{n}x$ crosses vertical lines 
$x=s$ at points of the form $\left(s,l+\frac{r}{n}\right)$, where $l$ is an integer
and $r\in\{0,1,2,...,n-1\}$. Numbers $m,n$ are relatively prime, therefore each value
of $r$ occurs precisely once for an appropriate value of $s\in\{0,1,2,...,n-1\}$.
Observe that the set of values $S_t(s)$, $s\in\{0,1,2,...,n-1\}$
is the same as the set of $t$-slopes at $0$ of the lines $y=\frac{m}{n}x+\frac{r}{n}$,
$r\in\{0,1,2,...,n-1\}$. Indeed an integer shift of coordinates allows us to view
each segment starting at $\left(s,\frac{m}{n}s\right)$ and ending at 
$\left(t+s,\frac{m}{n}(t+s)\right)$ as a segment starting at $\left(0,\frac{r}{n}\right)$
and ending at $\left(t,\frac{m}{n}t+\frac{r}{n}\right)$ for an appropriate value
of $r$. 
Line $y=\frac{m}{n}x$ crosses the vertical line $x=t$ at a point $l+\frac{u}{n}$
with $l$ an integer and $u\in\{1,2,...,n-1\}$. All lines $y=\frac{m}{n}x+\frac{r}{n}$
with $r\in\{0,1,2,...,n-u-1\}$ have the same $t$-slope at $0$ as the line $y=\frac{m}{n}x$.
The $t$-slope jumps up by $1$ for $r=n-u$ and keeps this value for all all the remaining 
$r\in\{n-u,n-u+1,...,n-1\}$. We conclude that function $S_t$ takes two values
$\left[\frac{m}{n}t+1\right]$ and $\left[\frac{m}{n}t+1\right]+1$. The first value 
is taken $n-u$ times and the second value $u$ times. We need to verify that
$n-u=n\left(\left[\frac{m}{n}t+1\right]-\frac{m}{n}t\right)$.
Indeed, $\frac{m}{n}t=l+\frac{u}{n}$, therefore $mt=nl+u$, $\left[\frac{m}{n}t+1\right]
=l+1$ and $n\left(\left[\frac{m}{n}t+1\right]-\frac{m}{n}t\right)=n(l+1)-mt=n-u$.

\noindent
{\bf Approximation to the identity by Fej\'er's kernel.} In the lemma that follows we quote a well known approximation to the identity property of the Fej\'er's kernel. We translate the original property from the group of one dimensional torus to the group of integers. For $f,g\in l^1(\Z)$ 
by $f*g$ we denote the convolution of $f$ and $g$, defined as
$$
f*g(k)=\sum_{l\in \Z} f(k-l)g(l).
$$ 
\begin{lem}\cite{Katz}\label{approximate_identity} 
If $K\in l^1(\Z)$, then
$$
\lim_{N\rightarrow \infty}\frac{1}{N}\left<K*\chi_{\left[0,N-1\right]},\chi_{\left[0,N-1\right]} \right>=\sum_{k\in \Z} K(k),
$$
where $\left< \cdot, \cdot \right>$ stands for the inner product of $l^2(\Z)$.
\end{lem}

\noindent
{\bf Invariance properties of the restricted kernel $K_R$.} Let us recall that kernel 
$$ 
K_R(\nu_1,\nu_2)=\frac{1}{R}\chi_{\Omega^c}(\nu_1)F(R(\nu_1-\nu_2))
\chi_{\Omega}(\nu_2),
$$
where $\Omega$ is a $\Z^2$ domain, $\nu_1,\nu_2\in \frac{1}{R}\Z^2$ and 
$F(\nu)=|\langle \mu(A)\phi,(\mu(A)\phi)_{S_{a,b}(\nu)} \rangle |^2$
was defined in (\ref{reproducing_kernel_z2}). This section presents 
invariance properties of kernel $K_R$ needed at various stages of the process of reduction. Let $\text{Aut}(\Z^2)=\text{Aut}(\frac{1}{R}\Z^2)=\{-1,+1\}^2
\rtimes S_2$, be the group of automorphisms of $\Z^2$ (or $\frac{1}{R}\Z^2$), i.e. the semi-direct product of sign changes and permutations of variables (see \cite{Mar} page 110). Let $\frac{1}{R}\Z^2 \rtimes \text{Aut}(\Z^2)$ be the group of affine transformations of $\frac{1}{R}\Z^2$ consisting of translations and 
 automorphisms of $\frac{1}{R}\Z^2$. The elements of $\frac{1}{R}\Z^2 \rtimes \text{Aut}(\Z^2)$ are represented as pairs $(\mu,h)$, where 
$\mu \in \frac{1}{R}\Z^2$, $h\in \text{Aut}(\Z^2)$. The group law
has the form $(\mu_1,h_1)(\mu_2,h_2)=(h_1(\mu_2)+\mu_1,h_1h_2)$.

\begin{lem}\label{invariance_properties}
Let $V,W\subset \R^2$ and let $\tau=(\mu,h)\in \frac{1}{R}\Z^2 \rtimes \text{Aut}(\Z^2)$. After the change of variables $\mu_1=\tau^{-1} (\nu_1)$, $\mu_2=\tau^{-1} (\nu_2)$ the restricted kernel $K_R$
$$
\frac{1}{R}\chi_V(\mu_1)\chi_{\Omega^c}(\mu_1)F(R(\mu_1-\mu_2))
\chi_{\Omega}(\mu_2)\chi_W(\mu_2)
$$
becomes
$$
\frac{1}{R}\chi_{\tau(V)}(\nu_1)\chi_{\tau(\Omega)^c}(\nu_1)
F(h^{-1}(R(\nu_1-\nu_2)))\chi_{\tau(\Omega)}(\nu_2)\chi_{\tau(W)}(\nu_2).
$$
For any $\Z^2$ domain $\Omega$ and any line segment $l\subset \partial \Omega$
it is possible to choose a transformation $\tau\in \frac{1}{R}\Z^2 \rtimes \text{Aut}(\Z^2)$ such that $\tau(l)$ is a segment of the graph of 
$H(x)=\frac{m}{n}x$, where $n\ge m$, $m\ge 0$, $n>0$, $m,n$ are relatively 
prime, and the portion of $\tau(\Omega)$ close to $\tau(l)$ is placed below $\tau(l)$, the portion of $\tau(\Omega)^c$ close to $\tau(l)$ is placed
above $\tau(l)$. If the line segment $l$ does not contain lattice points in its interior,
then we may assume that $\tau(l)=\{(x,H(x))\,|\,0\le x \le n\}$.
\end{lem}

\noindent
{\bf Proof.} The proof of the first part, the formula for the coordinate change is a straightforward computation which makes use of the fact that $\tau^{-1}=(-h^{-1}(\mu),h^{-1})$. The proof of the second part follows the process of inspecting the list of all possible positions of the segment $l$, which takes into account 
the placement of $\Omega$ and $\Omega^c$ relative to $l$, and indicating in each case the coordinate change $\tau$ needed to accomplish the target. It is also a direct 
computation.

\noindent
{\bf Asymptotic limits of the restricted kernel $K_R$.} The following two lemmas 
are the main technical tools behind the reduction process of general lattice domains to vertical strip domains. In the first lemma we deal with separated 
supports of variables $\nu_1,\nu_2$. The second lemma is more delicate,
it treats the case of variables $\nu_1,\nu_2$ restricted to bounded cones
located outside and inside $\Omega$, with their sides being segments of 
rational lines, and their common vertex being a lattice point of the boundary of $\Omega$. 
\noindent

\begin{lem}\label{separated_supports}
Suppose that there is $\delta>0$ such that for all sufficiently large $R$ 
nonnegative kernel $G_R(\nu_1,\nu_2),
\,\nu_1,\nu_2\in \frac{1}{R}\Z^2$ satisfies
$$
G_R(\nu_1,\nu_2)\le K_R(\nu_1,\nu_2)\text{, and }G_R(\nu_1,\nu_2)=0
\text{ for }|\nu_1-\nu_2|<\delta.
$$
Then
$$
\lim_{R\rightarrow \infty}\sum_{\nu_1,\nu_2\in \frac{1}{R}\Z^2}
G_R(\nu_1,\nu_2)=0.
$$
\end{lem}

\noindent
{\bf Proof.} For each $\nu_2 \in \frac{1}{R}\Z^2 \cap \Omega$ we have
\begin{align}
\sum_{\nu_1\in \frac{1}{R}\Z^2}G_R(\nu_1,\nu_2)&\le
\sum_{\substack{\nu_1\in \frac{1}{R}\Z^2 \\
|\nu_1-\nu_2|\ge \delta}}K_R(\nu_1,\nu_2)\le \frac{1}{R}
\sum_{\substack{\nu_1\in \frac{1}{R}\Z^2 \\
|\nu_1-\nu_2|\ge \delta}}\chi_{\Omega^c}(\nu_1)F(R(\nu_1-\nu_2))
\chi_\Omega (\nu_2) \nonumber \\
&\le \frac{\chi_\Omega (\nu_2)}{R}
\sum_{\substack{\nu\in \frac{1}{R}\Z^2 \\
|\nu|\ge \delta}}F(R\nu)=\frac{\chi_\Omega (\nu_2)}{R}
\sum_{\substack{\lambda\in \Z^2 \\
|\lambda|\ge \delta R}}F(\lambda)\le 
\frac{\chi_\Omega (\nu_2)}{\delta R^2}
\sum_{\substack{\lambda\in \Z^2 \\
|\lambda|\ge \delta R}}F(\lambda)|\lambda|.
\nonumber\end{align}
Let $Q_R=[\frac{-1}{2R},\frac{1}{2R}]\times [\frac{-1}{2R},\frac{1}{2R}]$. 
We observe that 
$$
\sum_{\nu_2\in \frac{1}{R}\Z^2}\frac{\chi_\Omega(\nu_2)}{R^2}=
\left| \bigcup_{\nu_2 \in \frac{1}{R}\Z^2 \cap \Omega}
\nu_2 + Q_R\right|
\le \left| \Omega^{\frac{1}{\sqrt{2}R}}\right|,$$
where for $r>0$ by $\Omega^r$ we denote $\{ u\in \R^2\,|\,\text{dist}(u,v)<r \text{ for some } v\in \Omega\}$.
Since $\sum_{\lambda\in \Z^2}F(\lambda)|\lambda|<\infty$ we obtain
$$
\sum_{\nu_1,\nu_2 \in \frac{1}{R}\Z^2}G_R(\nu_1,\nu_2)\le
\frac{1}{\delta}\left| \Omega^{\frac{1}{\sqrt{2}R}}\right|
\sum_{\substack{\lambda\in \Z^2 \\
|\lambda|\ge \delta R}}F(\lambda)|\lambda|
\ra_{R \rightarrow \infty} 0.
$$

\begin{lem}\label{supports_in_cones}
Suppose that there are two closed, bounded cones $C_1,C_2$, $C_1\subset 
\overline{\Omega^c}$, $C_2\subset \Omega$, with their sides being segments of rational lines, their apertures smaller than $\pi$ and their intersection being their common vertex, a lattice point of the boundary of $\Omega$. Suppose also that for all sufficiently large $R$ nonnegative kernel $G_R(\nu_1,\nu_2),
\,\nu_1,\nu_2\in \frac{1}{R}\Z^2$ satisfies
$$
G_R(\nu_1,\nu_2)\le K_R(\nu_1,\nu_2)\text{, and }G_R(\nu_1,\nu_2)=0
\text{ for }(\nu_1,\nu_2)\notin C_1\times C_2.
$$
Then
$$
\lim_{R\rightarrow \infty}\sum_{\nu_1,\nu_2\in \frac{1}{R}\Z^2}
G_R(\nu_1,\nu_2)=0.
$$
\end{lem}

\noindent
{\bf Proof.} We observe that in view of Lemma \ref{separated_supports} we may substitute cones $C_1,C_2$ by restricted cones $C_1^\delta=C_1\cap D_\delta (v)$,
$C_2^\delta=C_2\cap D_\delta (v)$, where $D_\delta (v)$ is the Euclidean disk 
with radius $\delta$ and center $v$, the common vertex of $C_1,C_2$. 
We may take the radius $\delta$ arbitrarily small. We choose 
a lattice line $l$ passing through $v$ and separating $C_1$ and $C_2$.
Lemma \ref{invariance_properties} allows us to assume that line $l$
has the form $H(x)=\frac{m}{n}x$, $m\ge 0$, $n>0$, $m\le n$,
$\text{gcd}(m,n)=1$, that the origin is the common vertex of $C_1,C_2$, 
and that $C_1$ lies above $l$ and $C_2$ below $l$.

If $m\ne 0$, i.e. line $l$ is not horizontal, then we extend the lattice $\Z^2$ to lattice $\tilde{\Lambda}$,
\begin{equation}
\tilde{\Lambda}=\left( \frac{1}{m}\Z\times \frac{1}{n}\Z \right)
\cup \left(\frac{1}{2m}, \frac{1}{2n}\right)+
\left( \frac{1}{m}\Z\times \frac{1}{n}\Z \right),
\label{lambda_tilde}\end{equation}
and kernel $F$ to kernel $\tilde{F}$,
\begin{equation}
\tilde{F}(\tilde{\lambda})=F(\lambda),
\label{F_tilde}\end{equation}
where $\lambda \in \Z^2$,
$Q=\left[ \frac{-1}{2},\frac{1}{2} \right)\times \left[ \frac{-1}{2},\frac{1}{2} \right)$, $\tilde{\lambda}\in \lambda + Q$. We observe that line $l$ becomes 
a horizontal line with respect to $\tilde{\Lambda}$. We obtain
\begin{align}
\frac{1}{R}
\sum_{\substack{\nu_1\in C_1^\delta \cap \frac{1}{R}\Z^2 \\
\nu_2\in C_2^\delta \cap \frac{1}{R}\Z^2}}F(R(\nu_1-\nu_2))
&\le\frac{1}{R}
\sum_{\substack{\tilde{\nu}_1\in C_1^\delta \cap \frac{1}{R}\tilde{\Lambda} \\
\tilde{\nu}_2\in C_2^\delta \cap \frac{1}{R}\tilde{\Lambda}}}
\tilde{F}(R(\tilde{\nu}_1-\tilde{\nu}_2))
=\frac{1}{R}
\sum_{\substack{\tilde{\lambda}_1\in C_1^{R\delta } \cap 
\tilde{\Lambda} \\
\tilde{\lambda}_2\in C_2^{R\delta }\cap \tilde{\Lambda}}}
\tilde{F}(\tilde{\lambda}_1-\tilde{\lambda}_2)
\nonumber\\
&\le \frac{1}{R}
\sum_{\tilde{\lambda}_1,\tilde{\lambda}_2\in \left( C_1\cup -C_2 \right)
^{R\delta } \cap 
\tilde{\Lambda}}
\tilde{F}(\tilde{\lambda}_1+\tilde{\lambda}_2)
\label{tilde_sum}
\end{align}
In the next step we bring the sum (\ref{tilde_sum}) to a computable form. We introduce an integer cone $0\le y\le \tilde{R}\delta $, 
$-My\le x \le My$, $x,y\in \Z$ with a sufficiently large 
aperture $M\in \Z$ in order to capture all points of 
$\left( C_1\cup -C_2 \right)^{R\delta } \cap \tilde{\Lambda}$
of (\ref{tilde_sum}). Constant $\tilde{\sigma}$ is the magnification 
factor needed to switch from $\tilde{\Lambda}$ to an integer lattice. 
Parameter $\tilde{R}=R\tilde{\sigma}$ of the integer cone accounts for this
magnification. We estimate (\ref{tilde_sum}) by 
\begin{equation}
\frac{\tilde{\sigma}}{\tilde{R}}\sum_{0\le y_1,y_2\le \tilde{R}\delta}
\sum_{\substack{-My_1\le x_1 \le My_1\\-My_2\le x_2 \le My_2}}
\breve{F}\left(\frac{1}{\tilde{\sigma}}(x_1+x_2,y_1+y_2)\right),
\label{integer_cone}\end{equation}
where $\breve{F}$ is $\tilde{F}$ adjusted to the new coordinate system
adapted to representing $l$ as a horizontal line. Let $x=x_1+x_2$, where $-My_1\le x_1 \le My_1$, $-My_2\le x_2 \le My_2$.
Then $-M(y_1+y_2)\le x \le M(y_1+y_2)$ and $x$ is represented as $x_1+x_2$
for at most $2M\text{min}\{ y_1,y_2\}+1$ pairs $x_1,x_2$. Therefore
we may estimate (\ref{integer_cone}) by
\begin{equation}
\frac{\tilde{\sigma}}{\tilde{R}}\sum_{0\le y_1,y_2\le \tilde{R}\delta}
\sum_{-M(y_1+y_2)\le x \le M(y_1+y_2)}
\breve{F}\left(\frac{1}{\tilde{\sigma}}(x,y_1+y_2)\right)
\left( 2M\text{min}\{ y_1,y_2\}+1\right).
\label{integer_cone_1}\end{equation}
Let $y=y_1+y_2$, where $0\le y_1,y_2\le \tilde{R}\delta$. Then 
$0\le y \le 2\tilde{R}\delta$ and $y$ is represented as $y_1+y_2$ for exactly 
$y+1$ pairs $y_1,y_2$. We may estimate (\ref{integer_cone_1}) by
\begin{equation}
\frac{\tilde{\sigma}}{\tilde{R}}\sum_{0\le y\le 2\tilde{R}\delta}
\sum_{-My\le x \le My}
\breve{F}\left(\frac{1}{\tilde{\sigma}}(x,y)\right)
\left( 2M y+1\right)(y+1).
\label{integer_cone_2}\end{equation}
Since $y+1\le 3\tilde{R}\delta$ for sufficiently large $\tilde{R}$, and
$2My+1\le 2M(y+1)$ we conclude that (\ref{integer_cone_2}) may be 
estimated by 
$$
\delta 6M\tilde{\sigma}\sum_{x,y\in \Z}\breve{F}\left(\frac{1}{\tilde{\sigma}}(x,y)\right)(y+1).
$$
This concludes the proof since constants $M,\tilde{\sigma}$ depend only on
cones $C_1,C_2$, condition $\Phi$ guaranties that $\sum_{x,y\in \Z}\breve{F}\left(\frac{1}{\tilde{\sigma}}(x,y)\right)(y+1)<\infty$, and we are allowed to take $\delta$ arbitrarily small.

\end{document}